\documentclass[12pt]{article}

\usepackage{amssymb}
\usepackage{amsthm}
\usepackage{latexsym}
\usepackage{epsfig}

\newcommand{\ben}{\begin{enumerate}}
\newcommand{\een}{\end{enumerate}}
\newcommand{\bde}{\begin{defn}}
\newcommand{\ede}{\end{defn}}
\newcommand{\bex}{\begin{exa}}
\newcommand{\eex}{\end{exa}}
\newcommand{\barr}{\begin{array}}
\newcommand{\earr}{\end{array}}
\newcommand{\btab}{\begin{tabular}}
\newcommand{\etab}{\end{tabular}}
\newcommand{\beq}{\begin{equation}}
\newcommand{\eeq}{\end{equation}}
\newcommand{\bea}{\begin{eqnarray*}}
\newcommand{\eea}{\end{eqnarray*}}
\newcommand{\bce}{\begin{center}}
\newcommand{\ece}{\end{center}}
\newcommand{\case}[4]{\left\{\barr{ll}#1&\mbox{#2}\\#3&\mbox{#4}\earr\right.}

\newcommand{\bib}{thebibliography}

\renewcommand{\qed}{\mbox{$\Diamond$}\vspace{\baselineskip}}
\newtheorem{theorem}{Theorem}[section]
\newtheorem*{higmansthm}{Higman's Theorem}

\newtheorem{proposition}[theorem]{Proposition}
\newtheorem{lemma}[theorem]{Lemma}

\newtheorem{corollary}[theorem]{Corollary}

\renewenvironment{proof}{\noindent {\bf Proof:}}{{\qed}}

\newcommand{\inv}{\mathop{\rm inv}}

\newcommand{\red}{\mathop{\rm red}}
\newcommand{\supp}{\mathop{\rm supp}}
\newcommand{\cl}{\mathop{\rm cl}}

\renewcommand{\Pr}{\mathop{\rm Prof}}
\newcommand{\Pa}{\mathop{\rm Part}}
\newcommand{\SPa}{\mathop{\rm SubPart}}

\newcommand{\Cong}{\equiv}
\newcommand{\Mod}{\mathop{\rm mod}\nolimits}

\newcommand{\mbar}{\overline{M}}
\newcommand{\pbar}{\overline{P}}
\newcommand{\ppbar}{\overline{P'}}
\newcommand{\ibar}{\overline{I}}

\newcommand{\ipbar}{\overline{I'}}

\newcommand{\A}{\mathcal{A}}
\newcommand{\zpm}{0/\mathord{\pm} 1}
\newcommand{\PrM}{(\nolinebreak[4]\Pr(M),\nolinebreak[4]\le\nolinebreak[4])}
\newcommand{\PaM}{(\nolinebreak[4]\Pa(M),\nolinebreak[4]\preceq\nolinebreak[4])}
\newcommand{\op}{\overline{P}}

\newcommand{\Gaa}{\put(0,0){\circle*{3}}}

\newcommand{\Gcb}{\put(20,10){\circle*{3}}}

\newcommand{\Gbc}{\put(10,20){\circle*{3}}}

\newcommand{\Gdd}{\put(30,30){\circle*{3}}}

\newcommand{\Gae}{\put(0,40){\circle*{3}}}
\newcommand{\Gbe}{\put(10,40){\circle*{3}}}
\newcommand{\Gce}{\put(20,40){\circle*{3}}}

\newcommand{\Gee}{\put(40,40){\circle*{3}}}
\newcommand{\Gfe}{\put(50,40){\circle*{3}}}
\newcommand{\Gge}{\put(60,40){\circle*{3}}}

\newcommand{\Gcf}{\put(20,50){\circle*{3}}}

\newcommand{\Gef}{\put(40,50){\circle*{3}}}

\newcommand{\Gag}{\put(0,60){\circle*{3}}}
\newcommand{\Gbg}{\put(10,60){\circle*{3}}}
\newcommand{\Gcg}{\put(20,60){\circle*{3}}}

\newcommand{\Geg}{\put(40,60){\circle*{3}}}
\newcommand{\Gfg}{\put(50,60){\circle*{3}}}
\newcommand{\Ggg}{\put(60,60){\circle*{3}}}
\newcommand{\Gddcf}{\put(30,30){\line(-1,2){10}}}
\newcommand{\Gddbe}{\put(30,30){\line(-2,1){20}}}

\newcommand{\Gcfag}{\put(20,50){\line(-2,1){20}}}
\newcommand{\Gbeag}{\put(10,40){\line(-1,2){10}}}
\newcommand{\Gddef}{\put(30,30){\line(1,2){10}}}
\newcommand{\Gddfe}{\put(30,30){\line(2,1){20}}}

\newcommand{\Gefgg}{\put(40,50){\line(2,1){20}}}
\newcommand{\Gfegg}{\put(50,40){\line(1,2){10}}}
\newcommand{\Gbcdd}{\put(10,20){\line(2,1){20}}}
\newcommand{\Gcbdd}{\put(20,10){\line(1,2){10}}}

\newcommand{\Gaabc}{\put(0,0){\line(1,2){10}}}
\newcommand{\Gaacb}{\put(0,0){\line(2,1){20}}}

\newcommand{\GbgLL}[2]{\Gbg \put(-3,64){\makebox(0,0){\hspace{#1}#2}}}
\newcommand{\GcgLL}[2]{\Gcg \put(7,64){\makebox(0,0){\hspace{#1}#2}}}

\newcommand{\GegLL}[2]{\Geg \put(27,64){\makebox(0,0){\hspace{#1}#2}}}
\newcommand{\GfgLL}[2]{\Gfg \put(37,64){\makebox(0,0){\hspace{#1}#2}}}
\newcommand{\GggLL}[2]{\Ggg \put(47,64){\makebox(0,0){\hspace{#1}#2}}}
\newcommand{\GaeLL}[2]{\Gae \put(-3,36){\makebox(0,0){\hspace{#1}#2}}}

\newcommand{\GceLL}[2]{\Gce \put(17,36){\makebox(0,0){\hspace{#1}#2}}}

\newcommand{\GeeLL}[2]{\Gee \put(37,36){\makebox(0,0){\hspace{#1}#2}}}

\newcommand{\GgeLL}[2]{\Gge \put(57,36){\makebox(0,0){\hspace{#1}#2}}}
\newcommand{\GddL}[2]{\Gdd \put(27,26){\makebox(0,0){\hspace{#1}#2}}}
\newcommand{\GbeL}[2]{\Gbe \put(-3,40){\makebox(0,0){\hspace{#1}#2}}}
\newcommand{\GfeL}[2]{\Gfe \put(47,36){\makebox(0,0){\hspace{#1}#2}}}
\newcommand{\GcfL}[1]{\Gcf \put(20,55){\makebox(0,0){#1}}}
\newcommand{\GefL}[1]{\Gef \put(40,55){\makebox(0,0){#1}}}
\newcommand{\GagL}[2]{\Gag \put(-13,64){\makebox(0,0){\hspace{#1}#2}}}
\newcommand{\GggL}[2]{\Ggg \put(57,56){\makebox(0,0){\hspace{#1}#2}}}
\newcommand{\GbcL}[2]{\Gbc \put(-3,20){\makebox(0,0){\hspace{#1}#2}}}
\newcommand{\GcbL}[2]{\Gcb \put(17,6){\makebox(0,0){\hspace{#1}#2}}}
\newcommand{\GaaL}[2]{\Gaa \put(-3,-4){\makebox(0,0){\hspace{#1}#2}}}
\newcommand{\Gaebg}{\put(0,40){\line(1,2){10}}}
\newcommand{\Gcebg}{\put(20,40){\line(-1,2){10}}}
\newcommand{\Geefg}{\put(40,40){\line(1,2){10}}}
\newcommand{\Ggefg}{\put(60,40){\line(-1,2){10}}}
\newcommand{\Gaefg}{\put(0,40){\line(5,2){50}}}
\newcommand{\Gcgeg}{\put(20,60){\line(1,0){20}}}
\newcommand{\Geggg}{\put(40,60){\line(1,0){20}}}
\newcommand{\Gcecg}{\put(20,40){\line(0,1){20}}}
\newcommand{\Geeeg}{\put(40,40){\line(0,1){20}}}
\newcommand{\Gceee}{\put(20,40){\line(1,0){20}}}
\newcommand{\Ggegg}{\put(60,40){\line(0,1){20}}}

\newcommand{\Geege}{\put(40,40){\line(1,0){20}}}

\makeatletter
\newfont{\footsc}{cmcsc10 at 8truept}
\newfont{\footbf}{cmbx10 at 8truept}
\newfont{\footrm}{cmr10 at 10truept}
\title{Profile classes and partial well-order for permutations}
\author{Maximillian M. Murphy\\
\small School of Mathematics and Statistics\\[-0.8ex]
\small University of St. Andrews\\[-0.8ex]
\small Scotland\\[-0.8ex]
\small \texttt{max@mcs.st-and.ac.uk}\\[1.6ex]
Vincent R. Vatter\thanks{Partially supported by an NSF VIGRE grant to the Rutgers University Department of Mathematics.}\\
\small Department of Mathematics\\[-0.8ex]
\small Rutgers University\\[-0.8ex]
\small USA\\[-0.8ex]
\small \texttt{vatter@math.rutgers.edu}}

\date{\small MR Subject Classifications: 06A06, 06A07, 68R15}

\begin{document}
\maketitle

\begin{abstract}
It is known that the set of permutations, under the pattern containment ordering, is not a partial well-order.  Characterizing the partially well-ordered closed sets (equivalently: down sets or ideals) in this poset remains a wide-open problem.  Given a $\zpm$ matrix $M$, we define a closed set of permutations called the profile class of $M$.  These sets are generalizations of sets considered by Atkinson, Murphy, and Ru\v{s}kuc.  We show that the profile class of $M$ is partially well-ordered if and only if a related graph is a forest.  Related to the antichains we construct to prove one of the directions of this result, we construct exotic fundamental antichains, which lack the periodicity exhibited by all previously known fundamental antichains of permutations.
\end{abstract}

\section{Introduction}

It is an old and oft rediscovered fact that there are infinite antichains of permutations with respect to the pattern containment ordering, so the set of all finite permutations is not partially well-ordered.  Numerous examples exist including Laver~\cite{laver}, Pratt~\cite{p:cpwdeqpsapq}, Tarjan~\cite{t:sunoqas}, and Speilman and B\'ona~\cite{sb:aiaop}.  In order to show that certain subsets of permutations are partially well-ordered, Atkinson, Murphy, and Ru\v{s}kuc~\cite{amr:pwocsop} introduced profile classes of $\zpm$ vectors (although they gave these classes a different name).  We extend their definition to $\zpm$ matrices, give a simple method of determining whether such a profile class is partially well-ordered, and add to the growing library of infinite antichains by producing antichains for those profile classes that are not partially well-ordered.  Finally, in Section 5 we generalize our antichain construction to produce exotic fundamental antichains.

The {\it reduction\/} of the length $k$ word $w$ of distinct integers is the $k$-permutation $\red(w)$ obtained by replacing the smallest element of $w$ by $1$, the second smallest element by $2$, and so on.  If $q\in S_k$, we write $|q|$ for the length $k$ of $q$ and we say that the permutation $p\in S_n$ {\it contains a $q$ pattern\/}, written $q\le p$, if and only if there is a subsequence $1\le i_1<\dots<i_k\le n$ so that $p(i_1)\dots p(i_k)$ reduces to $q$.  Otherwise we say that $p$ is $q$-avoiding and write $q\not\le p$.  The problem of enumerating $q$-avoiding $n$-permutations has received much attention recently, see Wilf~\cite{w:pip} for references.

The relation $\le$ is a partial order on permutations.  Recall that the partially ordered set $(X,\le)$ is said to be partially well-ordered if it contains neither an infinite properly decreasing sequence nor an infinite antichain (a set of pairwise incomparable elements).  Since $|q|<|p|$ whenever $q\le p$ with $q\neq p$, no set of permutations may contain an infinite properly decreasing sequence, so a set of permutations is partially well-ordered if and only if it does not contain an infinite antichain.

If $X$ is any set of permutations, we let $\A(X)$ denote the set of finite permutations that avoid every member of $X$.  We also let $\cl(X)$ denote the {\it closure\/} of $X$, that is, the set of all permutations $p$ such that there is a $q\in X$ that contains $p$.  We say that the set $X$ is {\it closed\/} (or that it is an {\it order ideal\/} or a {\it down-set\/}) if $\cl(X)=X$.  Now that we have the notation, we state another result from Atkinson et al.~\cite{amr:pwocsop}.

\begin{theorem}\label{singleperm}\cite{amr:pwocsop}
Let $p$ be a permutation.  Then $\A(p)$ is partially well-ordered if and only if $p\in\{1,12,21,132,213,231,312\}$.
\end{theorem}

We will rely heavily on the result of Higman~\cite{higman} that the set of finite words over a partially well-ordered set is partially well-ordered under the subsequence ordering.  More precisely, if $(X,\le)$ is a poset, we let $X^{\ast}$ denote the set of all finite words with letters from $X$.  Then we say that $a=a_1\dots a_k$ is a subsequence of $b=b_1\dots b_n$ (and write $a\le b$) if there is a subsequence $1\le i_1<\dots<i_k\le n$ such that $a_j\le b_{i_j}$ for all $j\in[k]$.

\begin{higmansthm}\cite{higman}
If $(X,\le)$ is partially well-ordered then so is $(X^{\ast},\le)$.
\end{higmansthm}

Actually, the theorem above is a special case of Higman's result, but it is all that we will need.

If $p\in S_m$ and $p'\in S_n$, we define the direct sum of $p$ and $p'$, $p\oplus p'$, to be the $(m+n)$-permutation given by
$$
(p\oplus p')(i)=\case{p(i)}{if $1\le i\le m$,}{p'(i-m)+m}{if $m+1\le i\le m+n$.}
$$
The skew sum of $p$ and $p'$, $p\ominus p'$, is defined by
$$
(p\ominus p')(i)=\case{p(i)+n}{if $1\le i\le m$,}{p'(i-m)}{if $m+1\le i\le m+n$.}
$$
Given a set $X$ of permutations, the {\it sum completion} of $X$ is the set of all permutations of the form $p_1\oplus p_2\oplus\dots\oplus p_k$ for some $p_1,p_2,\dots,p_k\in X$, and the {\it strong completion\/} of $X$ is set of all permutations that can be obtained from $X$ by a finite number of $\oplus$ and $\ominus$ operations.  The following result is given in \cite{amr:pwocsop}.

\begin{proposition}\label{completion}\cite{amr:pwocsop}
If $X$ is a partially well-ordered set of permutations, then so is the strong completion of $X$.
\end{proposition}

For example, this proposition shows that the set of layered permutations is partially well-ordered, as they are precisely the sum completion of the chain $\{1, 21, 321, \dots\}$.  Similarly, the set of separable permutations, the strong completion of the single permutation $1$, is partially well-ordered.

\section{Profile classes of $\zpm$ matrices}

This section is devoted to introducing the central object of our consideration: profile classes.  We begin with notation.  If $M$ is an $m\times n$ matrix and $(i,j)\in[m]\times[n]$, we denote by $M_{i,j}$ the entry of $M$ in row $i$ and column $j$.  For $I\subseteq[m]$ and $J\subseteq[n]$, we let $M_{I\times J}$ stand for the submatrix $(M_{i,j})_{i\in I,j\in J}$.  We write $M^t$ for the transpose of $M$.

Given a matrix of size $m\times n$, we define the its {\it support\/}, $\supp(M)$, to be the set of pairs $(i,j)$ such that $M_{i,j}\neq 0$.  The permutation matrix corresponding to $p\in S_n$, $M_p$, is then the $n\times n$ $0/1$ matrix with $\supp(M_p)=\{(i,p(i)) : i\in[n]\}$.

If $P$ and $Q$ are matrices of size $m\times n$ and $r\times s$ respectively, we say that $P$ contains a $Q$ pattern if there is a submatrix $P'$ of $P$ of the same size as $Q$ such that for all $(i,j)\in[r]\times[s]$,
$$
Q_{i,j}\neq 0\mbox{ implies }P'_{i,j}=Q_{i,j}.
$$

(Note that we have implicitly re-indexed the support of $P'$ here.)  We write $Q\le P$ when $P$ contains a $Q$ pattern and $Q\not\le P$ otherwise.  If $q$ and $p$ are permutations then $q\le p$ if and only if $M_q\le M_p$.  F\H{u}redi and Hajnal studied this ordering for $0/1$ matrices in \cite{fh:ds}.

We define the {\it reduction\/} of a matrix $M$ to be the matrix $\red(M)$ obtained from $M$ by removing the all-zero columns and rows.  Given a set of ordered pairs $X$ let $\Delta(X)$ denote the smallest $0/1$ matrix with $\supp(\Delta(X))=X$.  If we are also given a matrix $P$, let $\Delta^{(P)}(X)$ denote the matrix of the same size as $P$ with $\supp(\Delta^{(P)}(X))=X$, if such a matrix exists.  If $Q$ is a $0/1$ matrix satisfying $\red(Q)=Q$ (for instance if $Q$ is a permutation matrix) then $Q$ is contained in a $0/1$ matrix $P$ if and only if there is a set $X\subseteq \supp(P)$ with $\red(\Delta(X))=Q$.

We say that $M$ is a {\it quasi-permutation matrix\/} if there is a permutation matrix $M'$ that contains an $M$ pattern or, equivalently, if $\red(M)$ is a permutation matrix.  If $M$ is a quasi-permutation matrix and $\supp(M)=\{(i_1,j_1),\dots, (i_\ell,j_\ell)\}$ with $1\le i_1<\dots<i_\ell$, we say that $M$ is {\it increasing\/} if $1\le j_1<\dots<j_\ell$ and {\it decreasing\/} if $j_1>\dots>j_\ell\ge 1$.  Hence increasing quasi-permutation matrices reduce to permutation matrices of increasing permutations and decreasing quasi-permutation matrices reduce to permutation matrices of decreasing permutations.

In their investigation of partially well-ordered sets of permutations, Atkinson, Murphy, and Ru\v{s}kuc~\cite{amr:pwocsop} defined the ``generalized $W$s'' as follows.  Suppose $v=(v_1,\dots,v_s)$ is a $\pm 1$-vector and that $P$ is an $n\times n$ permutation matrix.  Then $P\in W(v)$ if and only if there are indices $1=i_1\le\dots\le i_{s+1}=n+1$ such that for all $\ell\in[s]$,
\ben
\item[(i)] if $v_\ell=1$ then $P_{[i_\ell,i_{\ell+1})\times[n]}$ is increasing,
\item[(ii)] if $v_\ell=-1$ then $P_{[i_\ell,i_{\ell+1})\times[n]}$ is decreasing.
\een
For example, the following matrix lies in $W\left(-1,1,1,-1\right)$ (the $0$ entries have been suppressed for readability).
$$
M_{532481697}=
\left(
\barr{rrrrrrrrr}
 & & & 		& 		1& & & 		& \\
 & &1& 		& 		& & & 		& \\
 &1& & 		& 		& & & 		& \\\hline
 & & & 		1& 		& & & 		& \\
 & & & 		& 		& & & 		1& \\\hline
1& & & 		& 		& & & 		& \\
 & & & 		& 		&1& & 		& \\\hline
 & & & 		& 		& & & 		&1 \\
 & & & 		& 		& &1& 		&
\earr
\right)
$$

Using Higman's Theorem, they obtained the following result.

\begin{theorem}\label{wk}\cite{amr:pwocsop}
For all $\pm 1$ vectors $v$, $(W(v),\le)$ is partially well-ordered.
\end{theorem}

Our goal in this section is to generalize the ``generalized $W$s'' and Theorem~\ref{wk}.  Suppose that $M$ is an $r\times s$ $\zpm$ matrix and $P$ is a quasi-permutation matrix.  An {\it $M$-partition of $P$\/} is a pair $(I,J)$ of multisets $I=\{1=i_1\le\dots\le i_{r+1}=n+1\}$ and $J=\{1=j_1\le\dots\le j_{s+1}=n+1\}$ such that for all $k\in[r]$ and $\ell\in[s]$,
\ben
\item[(i)] if $M_{k,\ell}=0$ then $P_{[i_k,i_{k+1})\times[j_\ell,j_{\ell+1})}=0$,
\item[(ii)] if $M_{k,\ell}=1$ then $P_{[i_k,i_{k+1})\times[j_\ell,j_{\ell+1})}$ is increasing,
\item[(iii)] if $M_{k,\ell}=-1$ then $P_{[i_k,i_{k+1})\times[j_\ell,j_{\ell+1})}$ is decreasing.
\een
For any $\zpm$ matrix $M$ we define the {\it profile class of $M$\/}, $\Pr(M)$, to be the set of all permutation matrices that admit an $M$-partition.  For instance, our previous example also lies in $
\Pr\left(
\barr{rrrr}
-1&-1&0&0\\
1&0&1&1
\earr
\right)$, as is illustrated below.
$$
M_{532481697}=\left(
\barr{rrr|rr|rr|rr}
 & & & 		& 		1& & & 		& \\
 & &1& 		& 		& & & 		& \\
 &1& & 		& 		& & & 		& \\
 & & & 		1& 		& & & 		& \\
\hline
 & & & 		& 		& & & 		1& \\
1& & & 		& 		& & & 		& \\
 & & & 		& 		&1& & 		& \\
 & & & 		& 		& & & 		&1 \\
 & & & 		& 		& &1& 		&
\earr
\right)
$$
Although we have arranged things so that profile classes are sets of permutation matrices, this will not stop us from saying that a permutation belongs to a profile class, and by this we mean that the corresponding permutation matrix belongs to the profile class.

Note that a matrix in $\Pr(M)$ may have many different $M$-partitions.  Also note that $W(v)=\Pr(v^t)$.  The profile classes of permutations defined by Atkinson~\cite{a:rp} fall into this framework as well: $p$ is in the profile class of $q$ if and only if $M_p\in\Pr(M_q)$.  (The wreath products studied in \cite{as:rpwp}, \cite{maximillian}, and briefly in the conclusion of this paper provide a different generalization of profile classes of permutations.)

Unlike the constructions they generalize, it is not true that the profile class of every $\zpm$ matrix is partially well-ordered.  For example, consider the Widderschin antichain $W=\{w_1,w_2,\dots\}$ given by
\bea
w_1&=&8,1\ |\ 5,3,6,7,9,4\ |\ |\ 10,11,2\\
w_2&=&12,1,10,3\ |\ 7,5,8,9,11,6\ |\ 13,4\ |\ 14,15,2\\
w_3&=&16,1,14,3,12,5\ |\ 9,7,10,11,13,8\ |\ 15,6,17,4\ |\ 18,19,2\\
&\vdots&\\
w_k&=&4k+4,1,4k+2,3,\dots,2k+6,2k-1\ |\\
&&2k+3,2k+1,2k+4,2k+5,2k+7,2k+2\ |\\
&&2k+9,2k,2k+11,2k-2,\dots,4k+5,4\ |\\
&&4k+6,4k+7,2
\eea
where the vertical bars indicate that $w_k$ consists of four different parts, of which the first part is the interleaving of $4k+4,4k+2,\dots,2k+6$ with $1,3,\dots,2k-1$, the second part consists of just six terms, the third part is the interleaving of $2k+9,2k+11,\dots,4k+5$ with $2k,2k-2,\dots,4$, and the fourth part has three terms.  Proofs that $W$ is an antichain may be found in \cite{amr:pwocsop,maximillian}, and this antichain is in fact a special case of our construction in Section 4, so Theorem~\ref{antichain} also provides a proof that $W$ forms an antichain.

Each $M_{w_k}$ has a $\left(
\barr{rr}
1&-1\\
-1&1
\earr
\right)$-partition: $(\{1,2k+3,4k+8\},\{1,2k+3,4k+8\})$.  For example,
$$
M_{w_2}=
\left(\barr{rrrrrr|rrrrrrrrr}
&&&&&&&&&&&1&&&\\
1&&&&&&&&&&&&&&\\
&&&&&&&&&1&&&&&\\
&&1&&&&&&&&&&&&\\
&&&&&&1&&&&&&&&\\
&&&&1&&&&&&&&&&\\\hline
&&&&&&&1&&&&&&&\\
&&&&&&&&1&&&&&&\\
&&&&&&&&&&1&&&&\\
&&&&&1&&&&&&&&&\\
&&&&&&&&&&&&1&&\\
&&&1&&&&&&&&&&&\\
&&&&&&&&&&&&&1&\\
&&&&&&&&&&&&&&1\\
&1&&&&&&&&&&&&&
\earr
\right)
\in\Pr\left(
\barr{rr}
1&-1\\
-1&1
\earr
\right).
$$
Therefore $\Pr\left(\barr{rr}1&-1\\-1&1\earr\right)$ is not partially well-ordered under the pattern containment ordering.

\thicklines
\setlength{\unitlength}{2pt}
\begin{figure}[t]
\begin{center}
\btab{c}
\begin{picture}(60,35)(0,30)
\GbgLL{30pt}{$x_1$} \GfgLL{30pt}{$x_2$}
\GaeLL{30pt}{$y_1$} \GceLL{30pt}{$y_2$} \GeeLL{30pt}{$y_3$} \GgeLL{30pt}{$y_4$}
\Gaebg\Gcebg
\Geefg\Ggefg
\Gaefg
\end{picture}
\\
Figure 1: $G\left(\barr{rrrr}1&1&0&0\\1&0&1&1\earr\right)$
\etab
\label{G-ex}
\end{center}
\end{figure}

If $M$ is an $r\times s$ $\zpm$ matrix we define the {\it bipartite graph\/} of $M$, $G(M)$, to be the graph with vertices $\{x_1,\dots,x_r\}\cup\{y_1,\dots,y_s\}$ and edges $\{(x_i,y_j):|M_{i,j}|=1\}$.  Figure~1 shows an example.  Our main theorem, proven in the next two sections, characterizes the matrices $M$ for which $\PrM$ is partially well-ordered in terms of the graphs $G(M)$:

\begin{theorem}\label{MAIN}
Let $M$ be a finite $\zpm$ matrix.  Then $\PrM$ is partially well-ordered if and only if $G(M)$ is a forest.
\end{theorem}

\section{When profile classes are partially well-ordered}

Our aim in this section is to prove the direction of Theorem~\ref{MAIN} that states that $\PrM$ is partially well-ordered if $G(M)$ is a forest.  In order to do this, we will need more notation.  In particular, we need to introduce two new sets of matrices, $\Pa(M)$ and $\SPa(M)$, and an ordering on them, $\preceq$.

We have previously defined $\Pr(M)$ to be the set of permutations matrices admitting an $M$-partition.  Now let $\Pa(M)$ consist of the triples $(P,I,J)$ where $P\in\Pr(M)$ and $(I,J)$ is an $M$-partition of $P$.  We let the other set, $\SPa(M)$, contain all triples $(P,I,J)$ where $P$ is a quasi-permutation matrix and $(I,J)$ is an $M$-partition of $P$.  Hence $\Pa(M)\subseteq\SPa(M)$.  

Suppose that $M$ is an $r\times s$ $\zpm$ matrix with $(P,I,J),(P',I',J')\in\SPa(M)$ where
\bea
I	&=&	\{i_1\le \dots\le i_{r+1}\},\\
J	&=&	\{j_1\le\dots\le j_{s+1}\},\\
I'	&=&	\{i_1'\le\dots\le i_{r+1}'\},\\
J'	&=&	\{j_1'\le\dots\le j_{s+1}'\}.
\eea
We write $(P',I',J')\preceq (P,I,J)$ if there is a set $X\subseteq\supp(P)$ such that $\red(\Delta(X))=\red(P')$ and for all $k\in[r]$ and $\ell\in[s]$,
$$
|X\cap ([i_k,i_{k+1})\times [j_\ell,j_{\ell+1}))|=|\supp(P')\cap ([i_k',i_{k+1}')\times [j_\ell',j_{\ell+1}'))|.
$$
Because $\Pa(M)\subseteq\SPa(M)$, we have also defined $\preceq$ on $\Pa(M)$.  It is routine to verify that $\preceq$ is a partial order on both of these sets.

The poset we are really interested in, $\PrM$, is a homomorphic image of $\PaM$.  Consequently, if for some $M$ we can show that $\PaM$ is partially well-ordered, then we may conclude that $\PrM$ is partially well-ordered.  This is similar to the approach Atkinson, Murphy, and Ru\v{s}kuc~\cite{amr:pwocsop} used to prove Theorem~\ref{wk}.  First we examine two symmetries of partition classes.

\begin{proposition}\label{transpose}
If $M$ is a $\zpm$ matrix then $(\Pa(M^t),\preceq)\cong\PaM$.
\end{proposition}
\begin{proof}
The isomorphism is given by $(P,I,J)\mapsto (P^t,J,I)$.
\end{proof}

Proposition~\ref{transpose} says almost nothing more than that for permutations $p$ and $q$, $q\le p$ if and only if $\inv(q)\le\inv(p)$, where here $\inv$ denotes the group-theoretic inverse.  Similarly, we could define the reverse of a matrix and see that $(\Pa(M),\preceq)\cong (\Pa(M'),\preceq)$ whenever $M$ and $M'$ lie in the same orbit under the dihedral group of order $4$ generated by these two operations.  In fact, we have the following more powerful symmetry.

\begin{proposition}\label{permute}
If $M$ and $M'$ are $\zpm$ matrices and $M'$ can be obtained by permuting the rows and columns of $M$ then $\PaM\cong(\Pa(M'),\preceq)$. 
\end{proposition}
\begin{proof}
By Proposition~\ref{transpose}, it suffices to prove this in the case where $M'$ can be obtained by permuting just the rows of $M$.  Furthermore, it suffices to show this claim in the case where $M'$ can be obtained from $M$ by interchanging two adjacent rows $k$ and $k+1$.  Let $(P,I=\{i_1\le\dots\le i_{r+1}\},J=\{j_1\le\dots\le j_{s+1}\})\in\Pa(M)$.  Define $P'$ by
\bea
P'_{[1,i_k)\times[n]}					&=&	P_{[1,i_k)\times[n]},\\
P'_{[i_k,i_k+i_{k+2}-i_{k+1})\times[n]}	&=&	P_{[i_{k+1},i_{k+2})\times[n]},\\
P'_{[i_k+i_{k+2}-i_{k+1},i_{k+2})\times[n]}	&=&	P_{[i_k,i_{k+1})\times[n]},\\
P'_{[i_{k+2},n]\times[n]}				&=&	P_{[i_{k+2},n]\times[n]},
\eea
and set
$$
I'=\{i_1\le\dots\le i_k\le i_k+i_{k+2}-i_{k+1}\le i_{k+2}\le\dots\le i_{r+1}\}.
$$
It is easy to check that $(P,I,J)\mapsto(P',I',J)$ is an isomorphism.
\end{proof}

The analogue of Proposition~\ref{transpose} for the poset $\PrM$ is true.  However, the analogue of Proposition~\ref{permute} fails in general.  For example, $\Pr\left(\barr{rrr}1&1&-1\earr\right)^t$ contains 21 permutations of length four, excluding only $3214$, $4213$, and $4312$, whereas $\Pr\left(\barr{rrr}1&-1&1\earr\right)^t$ is without $2143$, $3142$, $3241$, $4132$, and $4231$.  Propositions~\ref{transpose} and \ref{permute} suggest (although they fall short of proving) that whether or not $\PaM$ is partially well-ordered depends only on the isomorphism class of $G(M)$, this hint was the original motivation for our main result, Theorem~\ref{MAIN}.  We are now ready to prove one direction of this theorem.

\begin{theorem}\label{matrices}
Let $M$ be a $\zpm$ matrix.  If $G(M)$ is a forest then $\PaM$ is partially well-ordered.
\end{theorem}
\begin{proof}
Let $M$ be an $r\times s$ $\zpm$ matrix satisfying the hypotheses of the theorem.  By induction on $|\supp(M)|$ we will construct two maps, $\mu$ and $\nu$, such that if $(P,I,J)\in\SPa(M)$ then
$$
\nu(M;P,I,J)		=	\nu_1(M;P,I,J)\dots\nu_{|\supp(P)|}(M;P,I,J)\in([r]\times[s])^{|\supp(P)|},
$$
and
$$
\mu(M;P,I,J)		=	\mu_1(M;P,I,J)\dots\mu_{|\supp(P)|}(M;P,I,J)
$$
is a word containing each element of $\supp(P)$ precisely once, thus specifying an order for us to read through the nonzero entries of $P$.  The other map, $\nu$, will then record which section of $P$ each of these entries lie in.  This is formalized in the first of three claims we make about these maps below.

\ben
\item[(i)] If $\nu_t(M;P,I,J)=(a,b)$ then $\mu_t(M;P,I,J)\in[i_a,i_{a+1})\times[j_b,j_{b+1})$.
\item[(ii)] If $1\le a_1<\dots<a_b\le |\supp(P)|$ then
$$
\mu(M;\Delta^{(P)}(\{\mu_{a_1}(M;P,I,J),\dots,\mu_{a_b}(M;P,I,J)\}),I,J)
=
\mu_{a_1}(M;P,I,J)\dots\mu_{a_b}(M;P,I,J).
$$
\item[(iii)] If $(P',I',J')\in\SPa(M)$ with $\nu(M;P',I',J')=\nu(M;P,I,J)$ then
$$
\red(P')=\red(P).
$$
\een

First we show that this is enough to prove the theorem.  Higman's Theorem tells us that in any infinite set of words from $([r]\times[s])^\ast$ there are two that are comparable.  Hence in every infinite subset of $\Pa(M)$, there are elements $(P',I',J')$ and $(P,I,J)$ such that $\nu(M;P',I',J')\le\nu(M;P,I,J)$.  Hence there are indices $1\le a_1<\dots<a_b\le |\supp(P)|$ so that
$$
\nu(M;P',I',J')=\nu_{a_1}(M;P,I,J)\dots\nu_{a_b}(M;P,I,J).
$$
Now let $X=\{\mu_{a_1}(M;P,I,J),\dots,\mu_{a_b}(M;P,I,J)\}$.  Claim (ii) implies that
$$
\mu(M;\Delta^{(P)}(X),I,J)=\mu_{a_1}(M;P,I,J)\dots\mu_{a_b}(M;P,I,J),
$$
and thus by claim (i) we have
\bea
\nu(M;\Delta^{(P)}(X),I,J)
&=&
\nu_{a_1}(M;P,I,J)\dots\nu_{a_b}(M;P,I,J),\\
&=&
\nu(M;P',I',J').
\eea
Hence claim (iii) shows that
$$
\red(\Delta^{(P)}(X))=\red(P').
$$
This implies that $P'\le P$.  The other part of what we need to conclude that $(P',I',J')\preceq (P,I,J)$ comes directly from claim (i).  Therefore $\Pa(M)$ does not contain an infinite antichain, as desired.

We also need to say a few words about the symmetries of these matrices.  Suppose that we have constructed $\mu(M;P,I,J)$, and thus $\nu(M;P,I,J)$, for every $(P,I,J)\in\SPa(M)$.  We would like to claim that this shows how to construct $\mu(M^t;P,I,J)$ for every $(P,I,J)\in\SPa(M^t)$.

Let $(P,I,J)\in\SPa(M^t)$, so $(P^t,J,I)\in\SPa(M)$.  We define $\mu(M^t;P,I,J)$ in the natural way by
$$
\mu_t(M^t;P,I,J)=(b,a)\mbox{ if and only if }\mu_t(M;P^t,J,I)=(a,b).
$$
Claim (i) then shows us how to define $\nu(M^t;P,I,J)$.  Now suppose that $1\le a_1<\dots<a_b\le|\supp(P)|$ and let $X=\{\mu_{a_1}(M^t;P,I,J),\dots,\mu_{a_b}(M^t;P,I,J)\}$.  By definition,
$$
\mu_t(M;\Delta^{(P)}(X),I,J)=(b,a)
$$
for $t\in[b]$, where
$$
(a,b)=\mu_t(M^t;(\Delta^{(P)}(X))^t,J,I),
$$
and $(a,b)=\mu_{a_t}(M;P^t,J,I)$ by claim (ii) for $M$.  This shows that $\mu_t(M;\Delta^{(P)}(X),I,J)=\mu_{a_t}(M;P,I,J)$, proving claim (ii).  Claim (iii) is easier: if $(P',I',J'),(P,I,J)\in\SPa(M^t)$ have $\nu(M^t;P',I',J')=\nu(M^t;P,I,J)$ then $\nu(M;(P')^t,J',I')=\nu(M;P^t,J,I)$ so $\red((P')^t)=\red(P^t)$ and thus $\red(P')=\red(P)$.

We would also like to know how to construct $\mu(\mbar;P,I,J)$ if $\mbar$ is obtained by permuting the rows and columns of $M$.  By our work above, it suffices to show this when $\mbar$ can be obtained from $M$ by interchanging rows $k$ and $k+1$.  Let $(P,I=\{i_1\le\dots\le i_{r+1}\},J=\{j_1\le\dots\le j_{s+1}\})\in\SPa(\mbar)$ and define $\pbar$ by
\bea
\pbar_{[1,i_k)\times[n]}					&=&	P_{[1,i_k)\times[n]},\\
\pbar_{[i_k,i_k+i_{k+2}-i_{k+1})\times[n]}		&=&	P_{[i_{k+1},i_{k+2})\times[n]},\\
\pbar_{[i_k+i_{k+2}-i_{k+1},i_{k+2})\times[n]}	&=&	P_{[i_k,i_{k+1})\times[n]},\\
\pbar_{[i_{k+2},n]\times[n]}				&=&	P_{[i_{k+2},n]\times[n]},
\eea
and set
$$
\ibar=\{i_1\le\dots\le i_k\le i_k+i_{k+2}-i_{k+1}\le i_{k+2}\le\dots\le i_{r+1}\}.
$$
Note that $(\pbar,\ibar,J)\in\SPa(M)$, so we can construct $\mu(M;\pbar,\ibar,J)$.  Suppose that $\mu_t(M;\pbar,\ibar,J)=(a,b)$.  We construct $\mu(\mbar;P,I,J)$ by
$$
\mu_t(\mbar;P,I,J)=
\left\{
\barr{ll}
(a,b)&\mbox{if $(a,b)\notin [i_k,i_{k+2})\times[n]$,}\\
(a+i_{k+2}-i_{k+1},b)&\mbox{if $(a,b)\in[i_k,i_{k+1})\times[n]$,}\\
(a-(i_{k+1}-i_k),b)&\mbox{if $(a,b)\in[i_{k+1},i_{k+2})\times[n]$.}
\earr\right.
$$
As usual, claim (i) shows us how to construct $\nu(\mbar;P,I,J)$.  Checking claims (ii) and (iii) is similar to what we did for the transpose, so we omit it.

We are now ready to begin constructing $\mu$ and $\nu$.  If $M=0$, then the only members of $\SPa(M)$ are triples of the form $(P,I,J)$ where $P=0$.  In this event we set $\nu(M;P,I,J)$ and $\mu(M;P,I,J)$ to the empty word, and claims (i)--(iii) hold quite trivially.

Otherwise $G(M)$ has at least one edge, so it contains a leaf.  By our previous work, we may assume that $(r,s)\in\supp(M)$ and $(r,\ell)\notin\supp(M)$ for all $\ell<s$.  In other words, the last row of $M$ is identically $0$ except in the bottom-right corner, where it contains either a $1$ or $-1$.  Our construction of $\mu$ and $\nu$ will depend on the operations used to put $M$ into this form but this is of no consequence to us since we have shown that any definition of $\mu$ and $\nu$ that satisfies (i)--(iii) suffices to prove the theorem.

Let $\mbar=M_{[r-1]\times[s]}$.  Also, for any $(P,I=\{i_1\le\dots\le i_{r+1}\},J=\{j_1\le\dots\le j_{s+1}\})\in\Pa(M)$, let $\pbar=P_{[1,i_r)\times[1,j_{s+1})}$ and $\ibar=\{i_1\le\dots\le i_{r}\}$.  We have that $(\pbar,\ibar,J)\in\SPa(\mbar)$, and thus by induction we have maps
\bea
\nu(\mbar;\pbar,\ibar,J)
&=&
\nu_1(\mbar;\pbar,\ibar,J)\dots\nu_{|\supp(\pbar)|}(\mbar;\pbar,\ibar,J)\in([r-1]\times[s])^{|\supp(\pbar)|},\\
\mu(\mbar;\pbar,\ibar,J)
&=&
\mu_1(\mbar;\pbar,\ibar,J)\dots\mu_{|\supp(\pbar)|}(\mbar;\pbar,\ibar,J),
\eea
that satisfy (i), (ii), and (iii).

Now let us build another map, $\mu^{(0)}(M;P,I,J)$ by reading $P$ from left to right.  In other words,
$$
\mu^{(0)}(M;P,I,J)
=
\mu_1^{(0)}(M;P,I,J)\dots\mu_{|\supp(P)|}^{(0)}(M;P,I,J),
$$
where $\mu_a^{(0)}(M;P,I,J)$ is the element of $\supp(P)-\{\mu_1^{(0)}(M;P,I,J),\dots,\mu_{a-1}^{(0)}(M;P,I,J)\}$ with least second coordinate.

Clearly $\mu^{(0)}(M;P,I,J)$ contains each entry of $\supp(P)$ precisely once.  We will now form $\mu(M;P,I,J)$ by rearranging the entries of $\mu^{(0)}(M;P,I,J)$ that also lie in $\supp(\pbar)$ according to $\mu(\mbar;\pbar,\ibar,J)$.  More precisely, suppose that the elements of $\supp(\pbar)$ appear in positions $1\le a_1<\dots<a_{|\supp(\pbar)|}\le \supp(P)$ of $\mu^{(0)}(M;P,I,J)$.  Then let
$$
\mu(M;P,I,J)=\mu_1(M;P,I,J)\dots\mu_{|\supp(P)|}(M;P,I,J),
$$
where
$$
\mu_b(M;P,I,J)
=
\case
{\mu_c(\mbar;\pbar,\ibar,J)}{if $b=a_c$,}
{\mu_b^{(0)}(M;P,I,J)}{otherwise (which occurs when $\mu_b^{(0)}(M;P,I,J)\notin\supp(\pbar)$).}
$$

By claim (i), this also defines $\nu(M;P,I,J)$.  It remains to check that these maps have the desired properties.  If $z\in\supp(P)$, we will briefly use the notation $P-z$ to denote the matrix obtained from $P$ by changing the entry at $z$ to $0$.  To show (ii), it suffices to show that
\begin{eqnarray*}
\lefteqn{\nu(M; P-\mu_a(M;P,I,J),I,J)=}
\\ & &
\nu_1(M;P,I,J)\dots \nu_{a-1}(M;P,I,J)\nu_{a+1}(M;P,I,J)\dots \nu_{|\supp(P)|}(M;P,I,J).
\end{eqnarray*}
There are two cases to consider.  If $\mu_{a}(M;P,I,J)\in\supp(\pbar)$, then let $b$ be such that
$$
\mu_{a}(M;P,I,J)=\mu_b(\mbar;\pbar,\ibar,J),
$$
and let $c$ be such that
$$
\mu_a(M;P,I,J)=\mu_c^{(0)}(M;P,I,J).
$$
Clearly we have
\begin{eqnarray*}
\lefteqn{\mu^{(0)}(M;P-\mu_a(M;P,I,J),I,J)=}
\\ & &
\mu_1^{(0)}(M;P,I,J)\dots\mu_{c-1}^{(0)}(M;P,I,J)\mu_{c+1}^{(0)}(M;P,I,J)\dots\mu_{|\supp(P)|}^{(0)}(M;P,I,J),
\end{eqnarray*}
and by induction,
\begin{eqnarray*}
\lefteqn{\mu(\mbar;\overline{P-\mu_a(M;P,I,J)},\ibar,J)=}
\\ & &
\mu_1(\mbar;\pbar,\ibar,J)\dots\mu_{b-1}(\mbar;\pbar,\ibar,J)\mu_{b+1}(\mbar;\pbar,\ibar,J)\dots\mu_{|\supp(\pbar)|}(\mbar;\pbar,\ibar,J).
\end{eqnarray*}
This implies the claim.  The other case, where $\mu_a(M;P,I,J)\notin\supp(\pbar)$, is easier.

We now have only claim (iii) to show.  Suppose to the contrary that $(P',I',J'),(P,I,J)\in\SPa(M)$ satisfy $\nu(M;P',I',J')=\nu(M;P,I,J)$ but $\red(P')\neq\red(P)$, and choose $P'$ and $P$ with $|\supp(P)|=|\supp(P')|$ minimal subject to this constraint.  If $(r,s)$ occurs in neither of these words then we are done because $\ppbar=P'$, $\pbar=P$, and
$$
\nu(\mbar;\ppbar,\ipbar,J)=\nu(M;P',I',J')=\nu(M;P,I,J)=\nu(\mbar;\pbar,\ibar,J),
$$
so by induction on $|\supp(M)|$, $\red(P')=\red(P)$, a contradiction.

Otherwise $(r,s)$ occurs in both $\nu(M;P',I',J')$ and $\nu(M;P,I,J)$.  This is the only part of our proof that depends on the sign of $M_{r,s}$.  Since both cases are similar, we will show only the case where $M_{r,s}=1$.  Let $a$ be the position of the last occurance of $(r,s)$ in $\nu(M;P',I',J')$ and $\nu(M;P,I,J)$, so for all $b>a$, $\nu_b(M;P',I',J')=\nu_b(M;P,I,J)\neq (r,s)$.  By our assumptions on $M$ and the construction of $\mu$ and $\nu$, we know that of all elements in $\supp(P')$, $\mu_a(M;P',I',J')$ has the greatest first coordinate.  We also know the analogous fact for $\mu_a(M;P,I,J)$.

Furthermore, by claims (i) and (ii), we get
$$
\nu(M;P'-\mu_a(M;P',I',J'),I',J')=\nu(M;P-\mu_a(M;P,I,J),I,J),
$$
so by our choice of $P$ and $P'$, we have
$$
\red(P'-\mu_a(M;P',I',J'))=\red(P-\mu_a(M;P,I,J)).
$$

Due to our construction of $\mu$ and $\nu$ and our choice of $a$, each of the entries $\mu_1(M;P',I',J'),$ $\dots,$ $\mu_{a-1}(M;P',I',J')$ lies to the upper-left of $\mu_a(M;P',I',J')$, that is, they have lesser first and second coordinates.  In addition all of the other entries, $\mu_{a+1}(M;P',I',J'),$ $\dots,$ $\mu_{|\supp(P')|}(M;P',I',J')$, lie to the upper-right of $\mu_a(M;P',I',J')$.  Completely analogously, we have the same facts for $(P,I,J)$.  This is enough to conclude that $\red(P')=\red(P)$, a contradiction, proving the theorem.
\end{proof}

Theorem~\ref{matrices} and Proposition~\ref{completion} together imply the following corollary.

\begin{corollary}\label{matrices2}
If $M$ is a finite $\zpm$ matrix and $G(M)$ is a forest then the strong completion of $\PrM$ is partially well-ordered.
\end{corollary}

\section{When profile classes are not partially well-ordered}

We have half of Theorem~\ref{MAIN} left to prove, and its proof will occupy this section.  We would like to show that if $M$ is a $\zpm$ matrix for which $G(M)$ is not forest, i.e., it contains a cycle, then $\PrM$ contains an infinite antichain.  Our construction will generalize the Widderschin antichain introduced in the second section.

First an overview.  We will begin by constructing a chain
$$
\left(\barr{r}1\earr\right)=\overline{P}_1\le\overline{P}_2\le\dots
$$
of permutation matrices, each formed by inserting a new $1$ into the previous matrix in a specified manner.  Then from $\op_n$ we will form the $(n+2)\times(n+2)$ permutation matrix $P_n$ by expanding the ``first'' and ``last'' entries of $\op_n$ into appropriate $2\times 2$ matrices.  Finally, we will show that there is some constant $K$ depending only on $M$ for which each $P_n$ with $n\ge K$ has a unique $M$-partition, and from this it will follow that $\{P_n : n\ge K\}$ forms an antichain.

Before we begin, we need to make a technical observation.  If $M'\le M$ then $\Pr(M')\subseteq\Pr(M)$, so we will assume throughout this section that $G(M)$ is precisely a cycle.  This requirement is not strictly necessary, but it will simply the proofs greatly.

Now we are ready to construct $\op_n$, which will be an $n\times n$ permutation matrix containing $\op_{n-1}$.  To the nonzero entries of $\op_n$ we attach three pieces of information:
\ben
\item[(i)] a number; the entry we insert into $\op_{n-1}$ in order to form $\op_n$ will receive number $n$,
\item[(ii)] a yearn, which must be one of top-left, top-right, bottom-left, or bottom-right, and
\item[(iii)] a nonzero entry of $M$.
\een
We call the resulting object a {\it batch\/}, which will help us keep it separate from the entries of $M$.

When thinking about these three pieces of information, it might be best to think about starting with an empty matrix partitioned into blocks corresponding to the cells of $M$.  We will insert the batches in the order given by their number.  Each batch will be inserted into the block corresponding to the entry of $M$ given by (iii).  Within this block, each entry will be placed --- with some restrictions --- in the corner given by its yearn, so we might say colorfully that each batch yearns toward a corner of its block.  This implies that if the entry of $M$ corresponding to a batch is a $1$, then the yearn of that batch must be either top-left or bottom-right.  Otherwise the yearn must be top-right or bottom-left.  Finally, the entries that successive batches correspond to by (iii) will trace out the cycle in $G(M)$.

We have already stated that $\op_1=\left(\barr{r}1\earr\right)$, but we have not specified properties (ii) and (iii) of batch number $1$.  We can choose to correspond with the first batch any nonzero entry of $M$, but for the purpose of being as concise as possible, let us always take it to correspond to the left-most nonzero entry on the first row of $M$.  Such an entry exists because $G(M)$ has been assumed not to have isolated vertices.  Upon fixing this entry of $M$, we have two choices for the yearn of the first batch (although, up to symmetry, the two choices result in the same antichain, see Figure 2 for an example of this).  Let us always assume that the first batch yearns right-ward (either bottom-right if the corresponding entry of $M$ is a $1$ or top-right if it is a $-1$).

Having completed the definition of the first batch, we move on to the second.  Since $G(M)$ is precisely a cycle, there is a unique nonzero entry of $M$ on the same row as the entry that the first batch corresponds to.  We will choose this entry to correspond to the second batch.  (We have a choice to take the entry in the same row or the entry in the same column, but again it turns out that these two result in symmetric antichains.)  Finally, we specify that the second batch be top-yearning if the first batch was top-yearning and bottom-yearning if the first batch was bottom-yearning.  This, together with the sign of the corresponding entry of $M$, determines the yearn.

Before describing where to insert the second batch into $\op_1$ to form $\op_2$, let us define the other batches.  The $n$th batch will correspond to a nonzero cell of $M$ that shares either a row or column with the $n-1$st batch, but is not the same cell that either the $n-1$st batch or the $n-2$nd batch correspond to.  Such an entry exists because $G(M)$ is an even cycle (since $G(M)$ is bipartite for any $M$).  If the $n$th batch shares a row with the $n-1$st batch then the $n$th batch will have the same vertical yearning as the $n-1$st batch, that is, it will be top-yearning if the $n-1$st batch is top-yearning, and it will be bottom-yearning if the $n-1$st batch is bottom-yearning.  If the $n$th batch shares a column with the $n-1$st batch, then the two must share the same horizontal yearning.  Together with the sign of the corresponding entry of $M$, this determines the yearn of the $n$th batch.

Now suppose that we have $\op_{n-1}$ and want to insert the $n$th batch.  Suppose that this batch corresponds to the the cell $(i,j)\in\supp(M)$.  Then our first requirements are that the batch must be inserted
\ben
\item[(1)] below all batches corresponding to matrix entries $(x,y)$ with $x<i$,
\item[(2)] above all batches corresponding to matrix entries $(x,y)$ with $x>i$,
\item[(3)] to the right of all batches corresponding to matrix entries $(x,y)$ with $y<j$, and
\item[(4)] to the left of all batches corresponding to matrix entries $(x,y)$ with $y>j$.
\een
These four restrictions are enough to insure that the $n$th batch ends up in the desired ``block'' of $\op_n$.  Now we need to insure that it ends up in the correct position within this block.  To this end we place the $n$th batch as far towards its yearning as possible subject to (1)-(4) and one additional condition.  The $n$th and $n-1$st batches share either a column or a row, and due to this they must also share either their horizontal or vertical yearning, respectively.  The additional condition is simply that the $n$th batch must not overtake the $n-1$st batch in this yearning.

For example, suppose that the $n$th batch has top-left yearn and that the $n$th and $n-1$st batches share a row, so the $n-1$st batch also yearns to be high.  Then the $n$th batch must be placed below $n-1$st batch, but otherwise, subject to (1)-(4), as high and far to the left as possible.

Once we have constructed $\op_n$, we form $P_n$ by replacing the first and last batches by $\left(\barr{rr}1&0\\0&1\earr\right)$ if that batch corresponds to an $1$ in $M$ and by $\left(\barr{rr}0&1\\1&0\earr\right)$ if that batch corresponds to a $-1$ in $M$.

Before beginning the proof that the $P$ matrices form an antichain we do a small example, constructing $\op_1,\op_2,\dots,\op_6$ for the matrix
$$
M=\left(\barr{rr}1&-1\\-1&1\earr\right).
$$
Let us take the first batch to correspond to entry $(1,1)$ of $M$ and to have bottom-right yearn.  As for any $M$, we have
$$
\op_1=\left(\barr{r}1\earr\right).
$$
The second batch then corresponds to entry $(1,2)$ of $M$.  Since the first and second batches share a row, the second batch must be bottom-yearning, and thus its yearn must be bottom-left because $M_{1,2}=-1$.  To place the second batch into $\op_2$, we note that conditions (1)-(4) simply state that the second batch must be placed to the right of the first batch.  The other requirements insist that the second batch be placed above the first batch, so we end up with
$$
\op_2=\left(\barr{rr}&{\bf 1}\\1&\earr\right).
$$
(Here we have made the second batch bold and, as usual, suppressed the 0s.)

The third batch then corresponds to entry $(2,2)$ of $M$.  It must yearn leftward because it shares a column with the second batch, and since $M_{2,2}=1$, this means that its yearn must be top-left.  Conditions (1)-(4) imply that the third batch must be below the first and second batches and to the right of the first batch.  The other requirements imply that the third batch must be to the right of the second batch.  Hence we have
$$
\op_3=\left(\barr{rrr}&1&\\1&&\\&&{\bf 1}\earr\right).
$$

The fourth batch then has top-right yearn, and must lie below all the previous batches, to the right of the first batch, and to the left of the second and third batches, so
$$
\op_4=\left(\barr{rrrr}&&1&\\1&&&\\&&&1\\&{\bf 1}&&\earr\right).
$$

The fifth batch, like the first batch, corresponds to entry $(1,1)$ of $M$.  It has the same yearn as the first batch, bottom-right, and must be to the left of batches $2$, $3$, and $4$, above batches $3$ and $4$, but otherwise as far down and to the right as possible.  We then have
$$
\op_5=\left(\barr{rrrrr}&&&1&\\1&&&&\\&{\bf 1}&&&\\&&&&1\\&&1&&\earr\right).
$$

Finally, the sixth batch corresponds to entry $(1,2)$ of $M$ and has bottom-left yearn, like the second batch.  When this batch is inserted into $\op_5$ we get
$$
\op_6=\left(\barr{rrrrrr}&&&&1&\\1&&&&&\\&&&{\bf 1}&&\\&1&&&&\\&&&&&1\\&&1&&&\earr\right).
$$

To get $P_6$, we replace the first batch with the $2\times 2$ identity matrix and the last batch with the $2\times 2$ anti-identity matrix, resulting in
$$
P_6=\left(\barr{rrrrrrrr}&&&&&&1&\\1&&&&&&&\\&1&&&&&&\\&&&&&1&&\\&&&&1&&&\\&&1&&&&&\\&&&&&&&1\\&&&1&&&&\earr\right).
$$

\begin{figure}[t]
\begin{center}
\btab{cc}
\epsfig{file=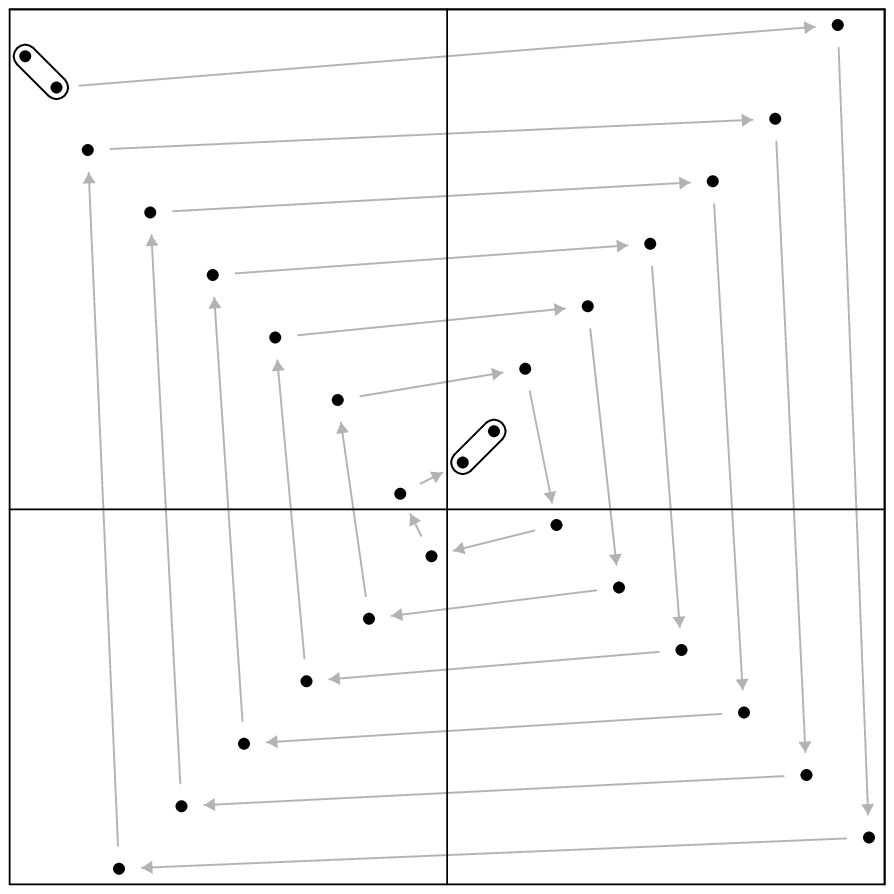, height=2.8in, width=2.8in}
&
\epsfig{file=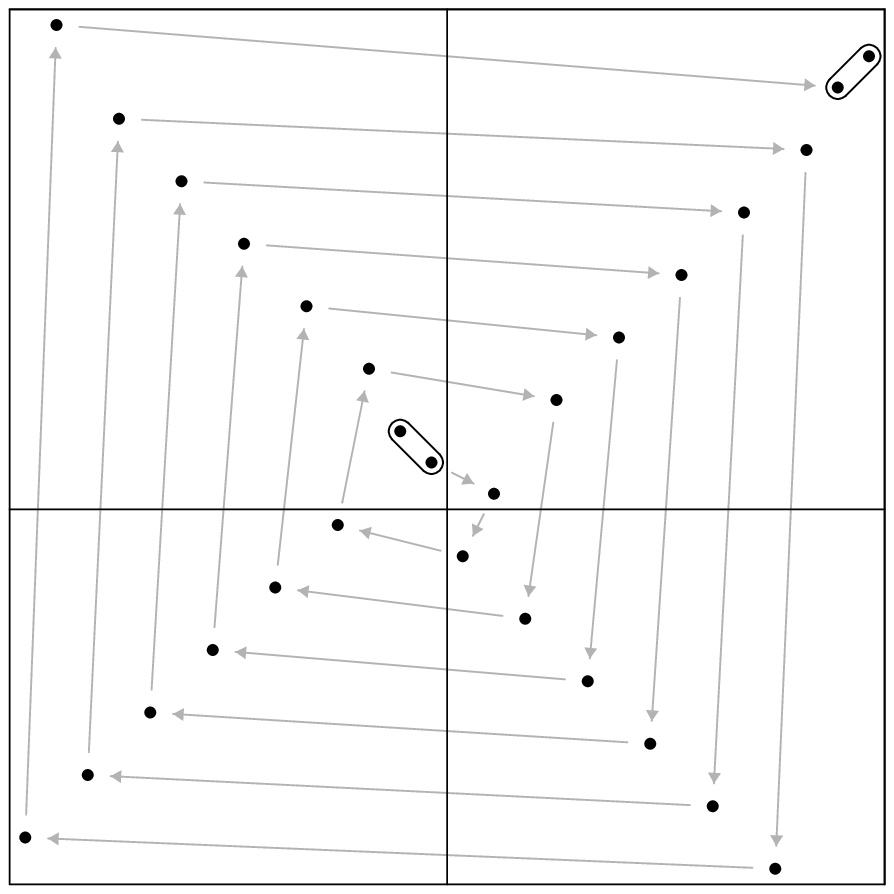, height=2.8in, width=2.8in}

\etab
\\
\end{center}
Figure 2: On the left we have a typical element of the Widderschin antichain initialized with bottom-right yearn, on the right it is initialized with top-left yearn.  Note that in this case the construction spins inward when initialized as on the left and outward when initialized as on the right, but the two resulting permutations are the same, up to symmetry.
\end{figure}

The matrix $P_{26}$ is shown on the left of Figure 2.  In the figure we have replaced 1s by dots and drawn an arrow from each batch to the subsequent batch.  It should be clear from that picture that we have constructed an antichain almost identical to the Widderschin antichain; the subset $\{P_9,P_{13},P_{17},P_{21},\dots\}$ is -- up to symmetry -- exactly the Widderschin antichain as we presented it in Section 2.

The matrix shown on the right of Figure 2 shows what we would get had we taken the yearn of the first batch to be top-left instead of bottom-right, and provides an example of our claim that the resulting matrices would be the same, up to symmetry.

Clearly the algorithm is well-defined since we have restricted ourselves to considering only cases where $G(M)$ is a cycle.  Almost as clearly, notice that successive batches cycle around $\supp(M)$.  Put more precisely, suppose that $G(M)$ is a cycle of length $c$.  Then the $m$th batch in $\op_n$ corresponds to the same entry of $M$ that the $(m+c)$th batch corresponds to. 

For the rest of our analysis we will restrict $M$ further, and assume that $M$ contains an even number of $-1$s.  That this can be done without loss is not completely obvious.  Suppose that $M$ has an odd number of $-1$s.  We form a new matrix $M'$ by replacing the $1$s in $M$ by $\left(\barr{rr}1&0\\0&1\earr\right)$, the $-1$s by $\left(\barr{rr}0&-1\\-1&0\earr\right)$, and the $0$s by the $2\times 2$ zero matrix.  It is easy to see that the profile classes of $M$ and $M'$ are identical, but we also need that $G(M')$ contains a unique cycle.

\begin{proposition}\label{evenminuses}
Let $M$ be a $\zpm$ matrix with an odd number of $-1$ entries, suppose that $G(M)$ is a cycle, and form $M'$ as described above.  Then $G(M')$ contains a unique cycle, twice the length of the cycle in $G(M)$.
\end{proposition}
\begin{proof}
First note that there is a natural homomorphism from $M'$ to $M$ that arises by identifying the $2\times 2$ blocks of $M'$ with the cells of $M$ that they came from.  This and the fact that every node of $G(M')$ has degree $2$ imply that $G(M')$ is either a cycle of twice the length of the cycle in $G(M)$ or two disjoint copies of the cycle in $G(M)$.  It is this latter case that we would like to rule out.

Now suppose that $M'$ is $r\times s$ and consider the $r\times s$ matrix $S$ given by $S_{i,j}=(-1)^{i+j}$.  For $r=s=4$, we have
$$
S=\left(\barr{rrrrr}
1&-1&1&-1\\
-1&1&-1&1\\
1&-1&1&-1\\
-1&1&-1&1
\earr\right).
$$

We can form $M'$ by changing entries of $S$ into $0$s.  A nice property of $S$ is that the length of any walk in $S$ in which diagonal steps are prohibited can be computed modulo $2$ by multiplying the entries of the start and end of the walk.  If this product is $1$ then then length of the walk is even, and if the product is $-1$ then the length of the walk is odd.  Clearly $M'$ also has this property for walks that begin and end at non-zero entries.

Now suppose that $G(M')$ contains two disjoint cycles and choose one of these cycles.  Clearly this cycle must contain one entry from each nonzero $2\times 2$ block of $M$, so the product of the matrix entries used in the cycle is $-1$, since $M$ had an odd number of $-1$s.  Now pair up the entries that lie on the same row of $M'$.  For any such pair, their product tells us whether the distance between them is odd or even.  Multiplying all such products together tells us whether the sum of all horizontal distances traversed by the cycle is odd or even.  Clearly, this result must be even.  However, this product will be the product of all entries used in the cycle, which we have assumed is $-1$.  Therefore this case is impossible, and the proposition is true.
\end{proof}

Hence we may assume that $M$ contains an even number of $-1$s.  This assumption will simplify our proofs, but it is worth noting that $M$ and $M'$ give rise to the same antichains; for example, see Figure 3.

\begin{figure}[t]
\begin{center}
\btab{cc}
\epsfig{file=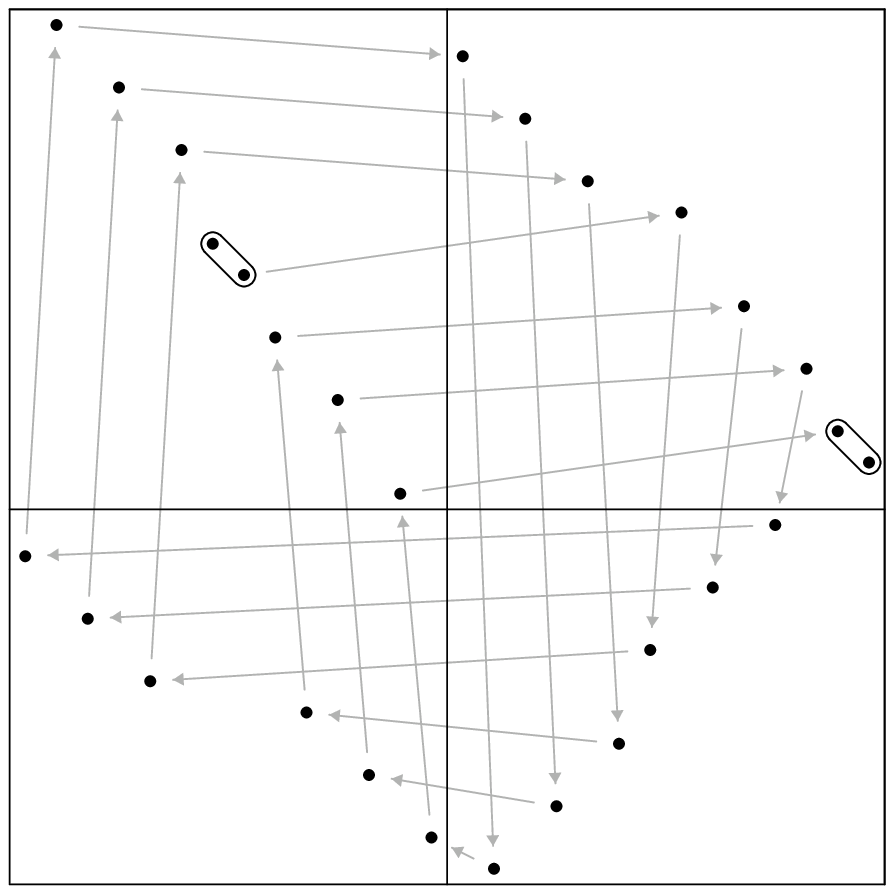, height=2.8in, width=2.8in}
&
\epsfig{file=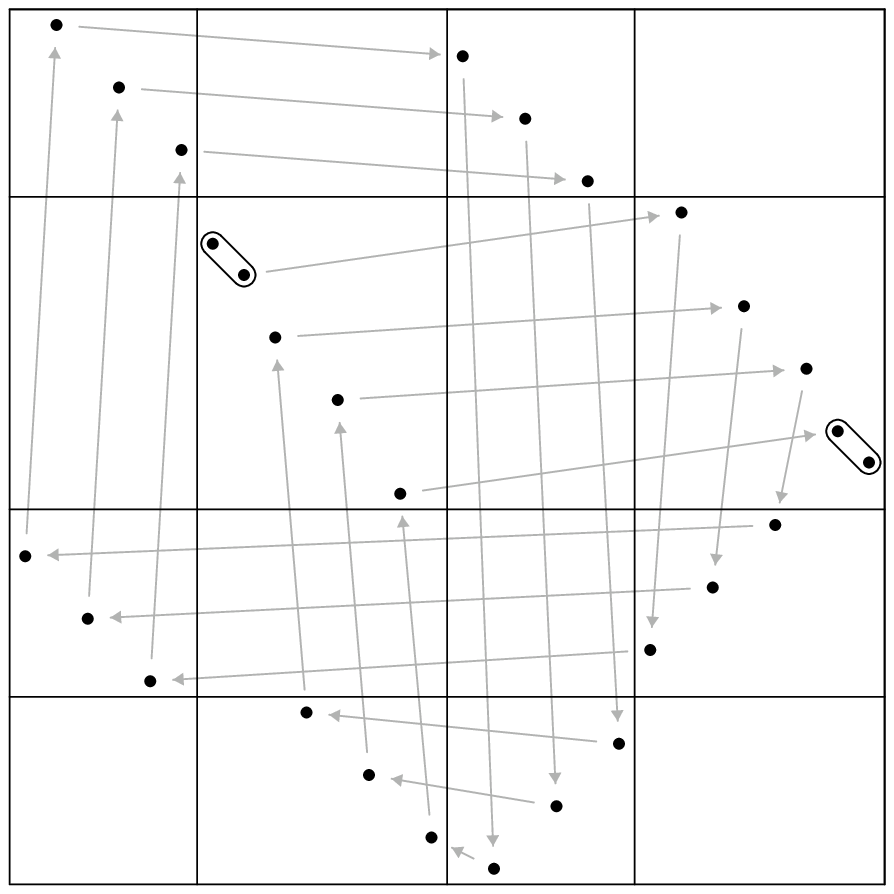, height=2.8in, width=2.8in}
\etab
\\
\end{center}
Figure 3: The permutation on the left is a permutation constructed by our algorithm to lie in $\Pr\left(\barr{rr}1&1\\1&-1\earr\right)$.  On the right, we have a permutation constructed by our algorithm to lie in the profile class of
$$
\left(\barr{rrrr}1&0&1&0\\0&1&0&1\\1&0&0&-1\\0&1&-1&0\earr\right).
$$
Unlike our previous examples, in this case the first batch was taken to correspond to the matrix entry $(2,2)$.  Notice that the two permutations are the same.
\end{figure}

Under this assumption we can conclude that any two batches that correspond to the same entry of $M$ share the same yearn.  Let us consider the vertical yearn only.  Note that it changes only when two successive batches correspond to cells of $M$ in the same column but with opposite signs.  Now if we sum over the columns of $M$ the number of $-1$ entries in each column we must get an even number.  Each time two successive batches correspond to cells of $M$ in the same column with the same sign, we either add $0$ or $2$ to our sum.  In the case that these batches correspond to cells of opposite sign, a $1$ is contributed.  Therefore the number of times that the vertical yearn changes during an entire cycle through $\supp(M)$ must be even.  The horizontal case is completely symmetric.

Because of this, one might think of set of batches that correspond to an entry of $M$ as ever more successfully progressing in the direction of their common yearn.

We now prove the major technical lemma we will need to establish that our algorithm produces antichains.

\begin{lemma}\label{unique-partition}
Suppose that $M$ is a $\zpm$ matrix with an even number of $-1$s and that $G(M)$ is a cycle of length $c$.  Then for $n\ge (c+1)c^2+1$, $P_n$ has a unique $M$-partition.
\end{lemma}
\begin{proof}
We begin by proving that under these hypotheses $\op_n$ has a unique $M$-partition, from which the claim for $P_n$ will follow rather easily.  First note that there is at least one $M$-partition of $\op_n$.  This $M$-partition comes from the correspondence between the batches of $\op_n$ and the entries of $M$.  Now suppose that we have another $M$-partition of $\op_n$.  We will use the verb ``allocate'' to differentiate this partition from the naturally arising partition just mentioned.  So, each batch of $\op_n$ corresponds to an entry of $M$ (this correspondence coming from the construction) and is allocated to an entry of $M$ (this allocation coming from the other $M$-partition we have been given).

Since $M$ has precisely $c$ nonzero entries and $n\ge (c+1)c^2+1$, we can find at least $c+2$ batches that correspond to the same cell of $M$ and are allocated to the same cell of $M$ (note that at this point we cannot assume that these two cells are the same).  We will call this set of batches isotope 1.  Because they are allocated to the same cell of $M$, the terms of isotope 1 form either an increasing or a decreasing sequence.  Suppose that the lowest numbered batch in isotope 1 is batch number $x$ and that the highest numbered batch is batch number $x+tc$ (so by our assumptions about the cardinality of isotope 1, $t\ge c+1$).  It should be clear from the principle that successive batches that correspond to the same entry of $M$ are increasingly successful in attaining their common yearn that isotope 1 also contains the batches $x,x+c,x+2c,\dots,x+tc$.

We know that all the batches in isotope 1 share the same yearn.  Without loss, we will assume that they are all top-yearning and that batch $x+1$ corresponds to an entry of $M$ on the same row as the entry that the batches of isotope 1 correspond to.  Then batch $x+1$ is lower than batch $x$, batch $x+c+1$ is higher than batch $x$ but lower than batch $x+c$, and in general batch $x+rc+1$ is lower than batch $x+rc$ but higher than batch $x+(r-1)c$ for all $r\in[t-1]$.  So, the batches numbered $x+rc+1$ for $r\in[t-1]$ horizontally separate the isotope 1 batches.  Therefore these batches must all be allocated to a nonzero entry of $M$ on the same row as the entry that the isotope 1 batches are allocated to.  However, these two isotopes must not be allocated to the same entry of $M$ on account of their non-monotonicity, and thus because $G(M)$ is a cycle, there is only one entry of $M$ that the batches numbered $x+rc+1$, $r\in[t-1]$, may be assigned to.  Let isotope 2 denote the set of all batches that correspond to and are allocated to the same entries of $M$ as these batches.

We now proceed in this manner, defining isotopes numbered 3 through $c+1$, each either vertically or horizontally separating the last.  In general isotope $i$ will contain at least $c+2-i$ batches.  Suppose that isotope $i-1$ contains all the batches numbered $x+rc+i-1$ where $r\in[s,t]$.  Then, excepting the $i=2$ case in which the existence of batch $x+tc+1$ is uncertain, isotope $i$ will contain the terms $x+rc+i-1$ where $r\in[s+1,t]$.  Thus we can be guaranteed that isotope $c+1$, the last isotope we will construct, is non-empty.

These isotopes must cycle around $M$, so the batches of isotope $c+1$ are allocated to the same entry of $M$ as the batches of isotope 1.  The isotopes must also contain a sequence of successive batches.  Furthermore, it is possible to determine the relative vertical and horizontal placement of the cells to which the isotopes are allocated, from which it follows that the batches in the isotopes are allocated to the cells to which they correspond.

It remains only to consider the batches that do not lie in the isotopes.  Note that if $c+1$ consecutive batches are allocated to the cells they correspond to, then also the batches immediately succeeding and preceding this sequence, if they exist, are allocated to the cells they correspond to.  For proof, suppose that such a sequence $y,y+1,\dots,y+c$ of batches is given and that, without loss, batches $y$ and $y+1$ are allocated to the same row and are both top-yearning.  Then batch $y+c+1$ lies vertically between batches $y$ and $y+c$, so it must be allocated to the same row as these batches (which is the same row that batch $y+1$ is allocated to).  Furthermore, it cannot horizontally separate batches $y$ and $y+c$, so it may not be allocated to the same column as these batches.  Since all rows of $M$ contain precisely two non-zero entries, this means that batch $y+c+1$ must be allocated to the same cell as batch $y+1$, and by our construction, this is the cell that batch $y+c+1$ corresponds to.

Therefore every batch in $\op_n$ must be allocated to the same cell that it corresponds to, and thus $\op_n$ has a unique $M$-partition.

Seeing that the same holds for $P_n$ is trivial.  Consider the first batch, which is expanded to form a $2\times 2$ matrix when we go from $\op_n$ to $P_n$.  Of the two nonzero entries in this $2\times 2$ matrix, one of them lies both horizontally and vertically between the other nonzero entry and batch $c+1$; call this entry interior, and then do the same for the last batch.  Removing the two interior entries gives $\op_n$ back, and we know that it has a unique $M$-partition.  But reinserting the interior entries cannot offer us any more possibilities for partitioning.
\end{proof}

\begin{figure}[t]
\begin{center}
\btab{cc}
\epsfig{file=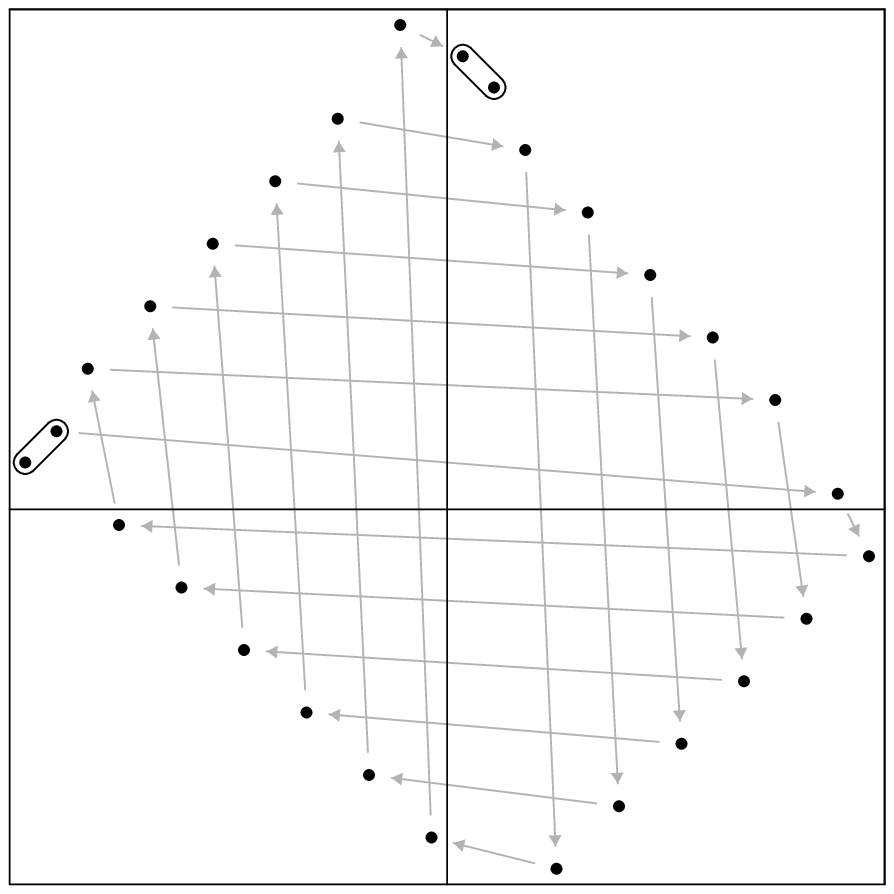, height=2.8in, width=2.8in}
&
\epsfig{file=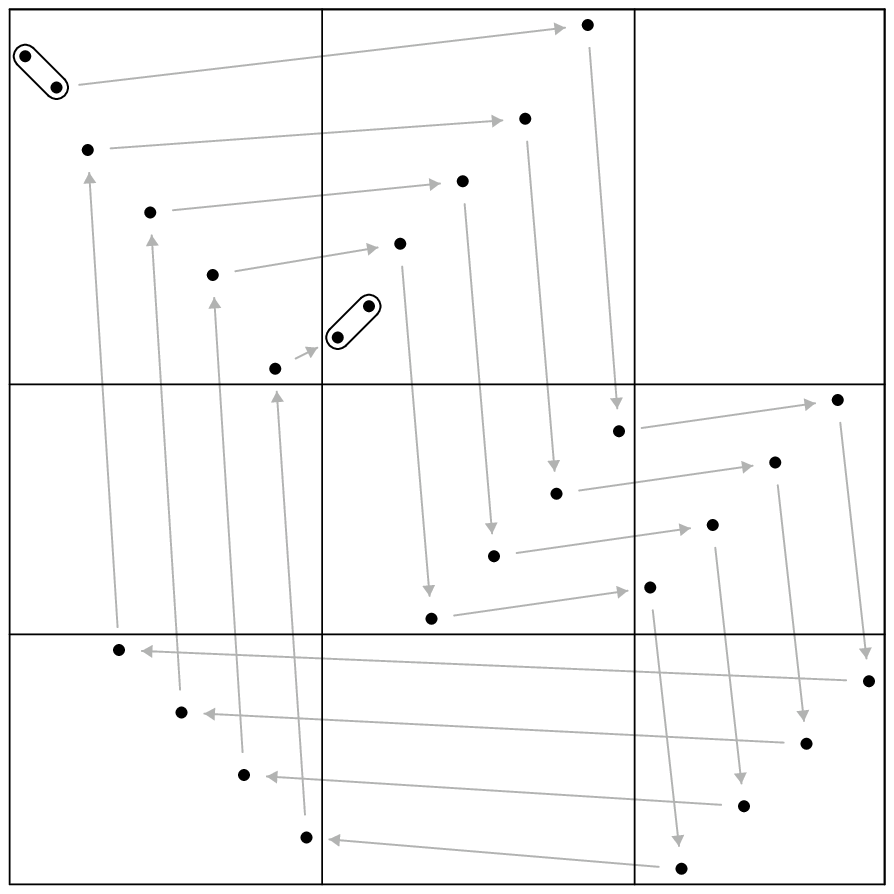, height=2.8in, width=2.8in}
\etab
\\
\end{center}
Figure 4: The permutation on the left is an element of the ``quasi-square antichain,'' introduced in \cite{maximillian} and readily constructed by our algorithm.  The permutation on the right comes from a matrix whose graph is a $6$-cycle.
\end{figure}

The main technical step now complete, we are ready to prove that the permutations we have constructed do indeed form antichains.

\begin{theorem}\label{antichain}
Let $M$ be a $\zpm$ matrix.  If $G(M)$ contains a cycle, then $\Pr(M)$ contains an infinite antichain, given by the permutations coming from $P_n$ for $n$ sufficiently large.
\end{theorem}
\begin{proof}
As we have already remarked, we may assume that $M$ contains an even number of $-1$s and that $G(M)$ is nothing but a cycle.  Let us assume this cycle is of length $c+1$, and that $n>m$ are both at least $c$ and large enough so that $P_m$ and $P_n$ have unique $M$-partitions (that we may make this assumption is the content of Lemma~\ref{unique-partition}).

We would like to show that $P_m$ and $P_n$ are incomparable.  Quite trivially, $P_n\not\le P_m$, because $P_n$ is larger than $P_m$, so it suffices to show that $P_m\not\le P_n$.  Suppose to the contrary that $P_m\le P_n$.  Then there is at least one submatrix of $P_n$ that reduces to $P_m$.  In this manner we get a one-to-one map from $\supp(P_m)$ into $\supp(P_n)$, or, as we will think of it, a map from the batches of $P_m$ into the batches of $P_n$ (with the first and last batches of $P_m$ possibly being mapped into more than one batch of $P_n$).  We begin by making two claims about this mapping:
\ben
\item[(i)] if batch $i$ of $P_m$ is mapped into batch $j$ of $P_n$ for some $i\in[2,m-1]$, then batch $i+1$ of $P_m$ is mapped into a batch of $P_n$ of numbered at most $j+1$, and
\item[(ii)] batch $1$ of $P_m$ must be mapped into batch $1$ of $P_n$.
\een

We begin with the proof of (i).  We may assume without loss that batches $i$ and $i+1$ of $P_m$ correspond to cells of $M$ that share a row, and that both batches are top-yearning.  First note since $P_m$ and $P_n$ have only one $M$-partition each, batch $i+1$ must be mapped to a batch of $P_n$ with number congruent to $j+1$ modulo $c$.  Furthermore, since batches $i$ and $i+1$ are both top-yearning and correspond to cells in the same row, batch $i+1$ lies below batch $i$, so it must be mapped to a batch below batch $j$.  These two restrictions leave only the possibilities we have allowed for.

Now we have to prove claim (ii).  Clearly batch $1$ of $P_m$ cannot be mapped into the last batch of $P_n$, so if the claim does not hold then this batch is mapped into two batches of $P_n$, say $r+1$ and $r+sc+1$ (by the uniqueness of $M$-partitions, these two batches must be congruent to $1$ modulo $c$).  Let us suppose that the first batch of $P_m$ has top-right yearn, and that the first and second batches of $P_m$ correspond to cells of $M$ that share a row.  Now consider where batch $2$ may be mapped to.  By the same argument we used in (i), batch $2$ must be mapped to a batch this lies below the batch $r+1$, so it must be mapped to a batch with number at most $r+2$.  We now follow the implication of (i) all the way around the cycle, to see that batch $c+1$ of $P_m$ must be mapped into a batch with number at most $r+c+1$.  This is a contradiction because either $r+c+1=r+sc+1$ (our mapping was supposed to be one-to-one) or $r+c+1<r+sc+1$ and thus batch $c+1$ is mapped to a batch that horizontally separates the two entries that batch $1$ was mapped to.

Having established (i) and (ii) we are almost done.  The first batch of $P_m$ must be mapped to the first batch of $P_n$, so the second batch of $P_m$ must be mapped to the second batch of $P_n$, and so on, until we conclude that the $m-1$st batch of $P_m$ must be mapped to the $m-1$st batch of $P_n$.  Now we have no options for the last batch.  Suppose without loss that the $m-1$st and $m$th batches of $P_m$ correspond to row-sharing cells of $M$, and that the $m-1$st batch is top-yearning.  Then the $m$th batch (which consists of two non-zero entries) must lie entirely below the $m-1$st batch.  This means that the $m$th batch of $P_m$ must be mapped to a batch of number at most $m$ in $P_n$.  Additionally, of course, the $m$th batch of $P_m$ may not be mapped into a batch that any other batch of $P_m$ has been mapped into, so we have reached a contradiction, proving the theorem.
\end{proof}

\section{Exotic Fundamental Antichains}

If we have an antichain of permutations $A$, then we may form infinitely more antichains from it by direct sums (or skew sums, or in several other ways).  For instance, $\{1324\oplus a : a\in A\}$ must also be an antichain.  But $\{1324\oplus a : a\in A\}$ is, at least intuitively, less interesting that $A$.  

In order to make this intuition precise, we say that an antichain $A$ is {\it fundamental\/} if its closure contains no antichains of the same size as $A$, except those that are subsets of $A$ itself.

Clearly $\{1324\oplus a : a\in A\}$ is not fundamental.  Note that some researchers (for example, Cherlin and Latka~\cite{cl:minimal} and Gustedt~\cite{gustedt}) call such antichains ``minimal.''

While we have no use for it, we would be remiss if we did not make note of the following result.  Surely the proof (or some generalization of it) has appeared in more than the two sources we cite.

\begin{proposition}\cite{gustedt, maximillian}
Let $X$ be a closed set of permutations.  If $X$ contains an infinite antichain, then it also contains an infinite fundamental antichain.
\end{proposition}

It can be shown that our construction from the previous section produces fundamental antichains.  However, this is a rather subtle point.  Consider our construction applied to the matrix
$$
M=\left(\barr{rr}1&-1\\-1&1\earr\right),
$$
and suppose we take the first batch to correspond to entry $(1,1)$ of $M$ and to have bottom-right yearn.  Our construction will then produce a sequence $P_1,P_2,\dots$ of permutation matrices.  Theorem~\ref{antichain} shows that the subset $\{P_n : n\ge 81\}$ forms an antichain, and indeed it is easy to check that $\{P_n : n\ge 9\}$ forms an antichain.

As we remarked in that section, the subset $\{P_n : 9\le n\Cong 1\ (\Mod 4)\}$ is, up to symmetry, exactly the Widderschin antichain as we presented it in Section 2, and this antichain is fundamental.  However, if we extend this antichain by adding $P_{10}$, it is no longer fundamental.

For proof of this, consider the permutation matrix $P_{10}'$ obtained from $P_{10}$ by removing one of the two entries coming from the first batch, shown below with its unique $M$-partition:
$$
P_{10}'=\left(\barr{rrrrr|rrrrrr}
&&&&&&&&&1&\\
1&&&&&&&&&&\\
&&&&&&&1&&&\\
&1&&&&&&&&&\\
&&&&&&1&&&&\\
&&&&&1&&&&&\\
&&&1&&&&&&&\\\hline
&&&&&&&&1&&\\
&&&&1&&&&&&\\
&&&&&&&&&&1\\
&&1&&&&&&&&
\earr\right).
$$
It is not hard to check that $P_{10}'\not\le P_n$ for all $13\le n\Cong 1\ (\Mod 4)$.  Therefore
$$
\{P_{10}'\}\cup\{P_n : 13\le n\Cong 1\ (\Mod 4)\}
$$
forms an antichain, which lies in the closure of
$$
\{P_{10}\}\cup\{P_n : 13\le n\Cong 1\ (\Mod 4)\}
$$
but is not a subset of this latter antichain.  Therefore this antichain is not fundamental, and in particular, $\{P_n : n\ge 9\}$ is not fundamental.

In general, suppose that $M$ is a $\zpm$ matrix for which $G(M)$ is precisely a cycle of length $c$.  If we fix some integer $d\in[c]$, the set
$$
\{P_n : n\mbox{ is sufficiently large and }n\Cong d\ (\Mod c)\}
$$
can be shown to form a fundamental antichain.

Up to this point, all fundamental antichains in the literature and all antichains produced by our algorithm as we have described it are periodic in some sense.  We aim in this section to convince the reader that our construction from the last section can be generalized to construct exotic fundamental antichains without this periodicity.

In the description of the construction we assumed that $G(M)$ was a cycle.  Suppose instead that we let $G(M)$ contain two or more cycles intersecting at a single vertex.  For example, let us take
$$
M=
\left(\barr{rrrrrr}
-1&1&-1&1&-1&1\\
1&-1&0&0&0&0\\
0&0&1&-1&0&0\\
0&0&0&0&1&-1
\earr\right),
$$
so $G(M)$ will be the graph depicted in Figure 5.

\thicklines
\setlength{\unitlength}{2pt}
\begin{figure}[t]
\begin{center}
\btab{c}
\begin{picture}(72,80)(-13,-14)
\GggL{30pt}{$x_4$}
\GfeL{30pt}{$y_6$}
\GddL{30pt}{$x_1$}
\GcbL{30pt}{$y_1$}
\GaaL{30pt}{$x_2$}
\GcfL{$y_4$}
\GefL{$y_5$}
\GbcL{30pt}{$y_2$}
\GbeL{30pt}{$y_3$}
\GagL{30pt}{$x_3$}
\Gddcf\Gddbe\Gcfag\Gbeag
\Gddef\Gddfe\Gefgg\Gfegg
\Gbcdd\Gcbdd
\Gaabc\Gaacb
\end{picture}
\\
Figure 5: $G\left(\barr{rrrrrr}
-1&1&-1&1&-1&1\\
1&-1&0&0&0&0\\
0&0&1&-1&0&0\\
0&0&0&0&1&-1
\earr\right)$
\etab
\end{center}
\end{figure}

Let us take the first batch to correspond to cell $(1,1)$ and to have bottom-right yearn.  Immediately we are faced with a predicament: there are five other entries on the first row of $M$, so which should we choose for the second batch?  This situation can be rectified by supplying additional input to the algorithm.  Let us supply a word $w$ on the letters $0$, $1$, and $2$ where $0$ means that the next several batches should trace out the left cycle, $(1,1)$, $(1,2)$, $(2,2)$, and $(2,1)$, $1$ means that the next several batches should trace out the middle cycle, $(1,3)$, $(1,4)$, $(3,4)$, and $(3,3)$, and $2$ means that the next several batches should trace out the cycle in columns $5$ and $6$.  Then our construction is once again well-defined, and a slight adaptation of the proofs in the last section would show that it still produces antichains.

To produce an aperiodic antichain we need only select an aperiodic word as $w$.  We define the (infinite) binary Thue-Morse word, $t$, by $t=\lim_{n\rightarrow\infty} u_n$ where $u_0=a$, $v_0=b$, and for $n\ge 1$, $u_n=u_{n-1}v_{n-1}$ and $v_n=v_{n-1}u_{n-1}$.  For example,
$$
u_6=abbabaabbaababbabaababbaabbabaab.
$$
Now we replace each occurrence of $abb$ by $2$, each occurrence of $ab$ by $1$ (after replacing the occurrences of $ab$), and each remaining occurance of $a$ by $0$ to get the word $w$ on the letters $0,1,2$.  Applying these substitutions to $u_6$, we obtain the word
$$
21020121012021.
$$
It is known (see, for example, Lothaire~\cite{lothaire1}) that the word $w$ is square-free.  An element of an antichain produced in this manner is shown in Figure 6.

To get an aperiodic fundamental antichain from this construction, we need only make sure to take elements $P_n$ for which the last batch always corresponds to the same cell of $M$.  For example, if we take all permutation matrices $P_n$ (for $n$ sufficiently large) produced by the operation described above, they will form an antichain, but not a fundamental one.  Instead if we take all elements $P_n$ where $n$ is sufficiently large and the last batch corresponds to entry $(1,1)$, this will be a fundamental antichain. 

\begin{figure}
\begin{center}
\epsfig{file=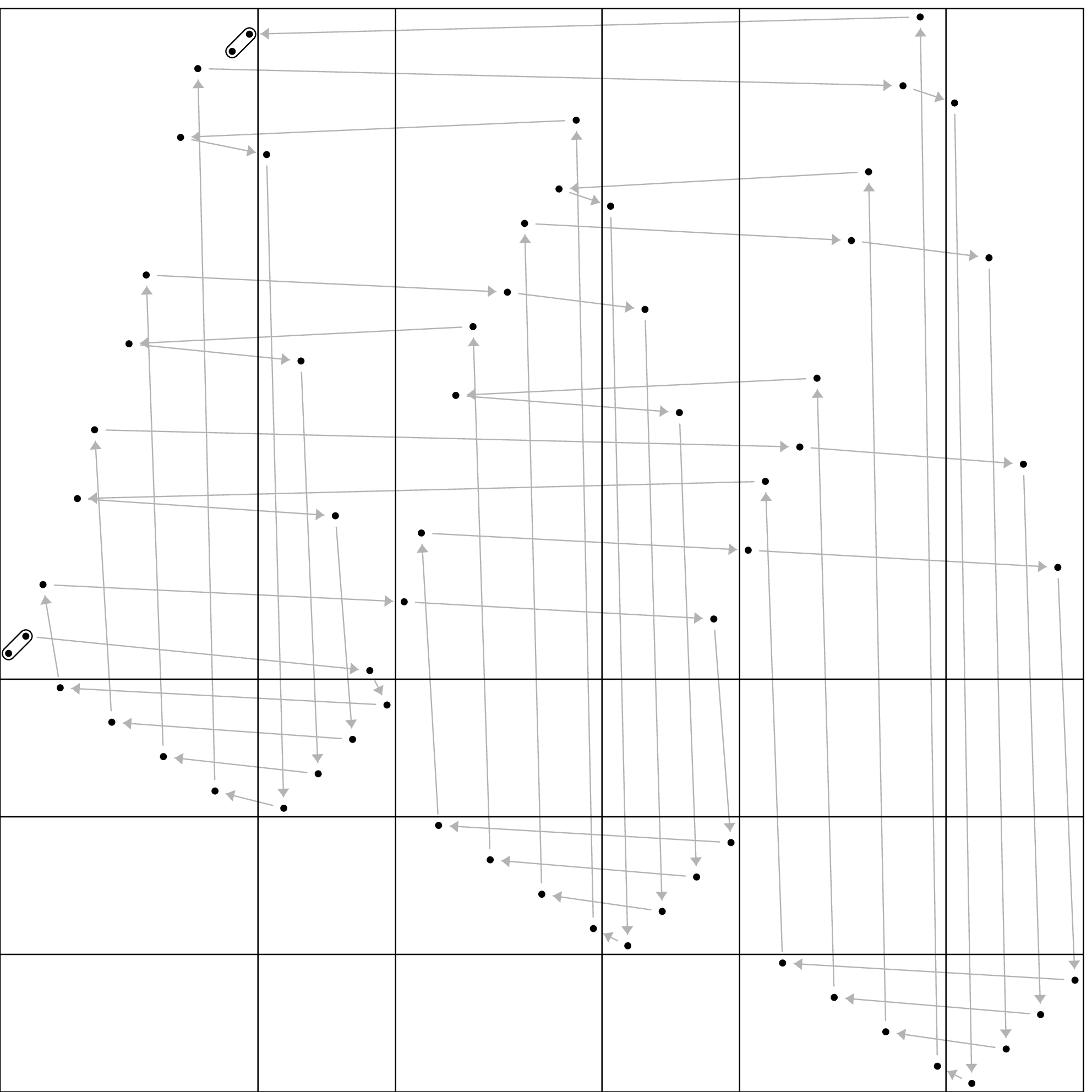, height=6.3in, width=6.3in}
\end{center}
Figure 6: An element of an aperiodic infinite antichain constructed from the Thue-Morse word.
\end{figure}

\begin{figure}
\begin{center}
\epsfig{file=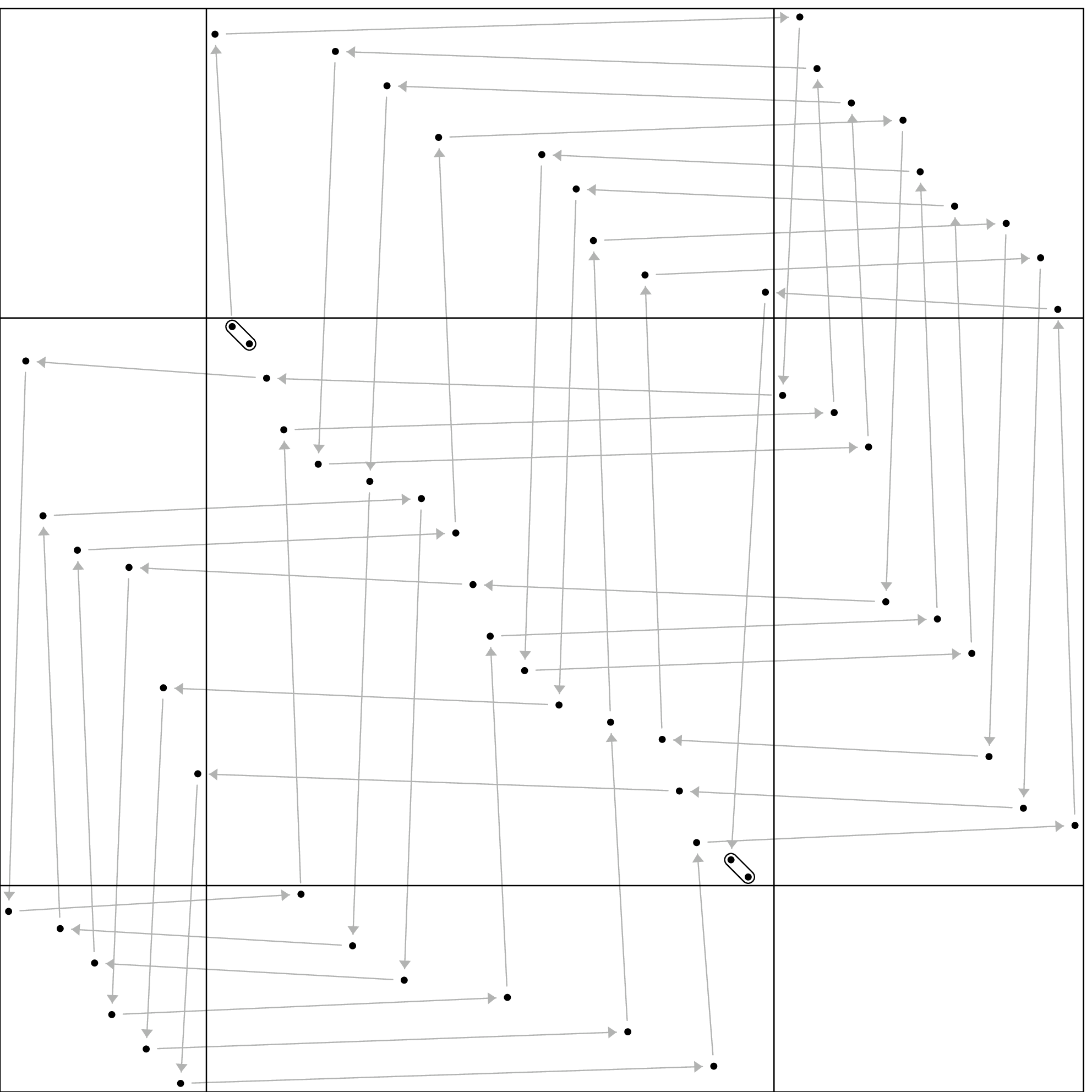, height=6.3in, width=6.3in}
\end{center}
Figure 7: An element of an aperiodic infinite antichain constructed from a matrix with two cycles and the Thue-Morse word.
\end{figure}

Even though $w$ contains three letters, with a little care we can use it to build an antichain in the profile class of matrix whose graph has only two cycles, for instance, the matrix
$$
M=\left(\barr{rrr}
0&1&1\\
1&1&1\\
1&1&0
\earr\right).
$$

\thicklines
\setlength{\unitlength}{2pt}
\begin{figure}[t]
\begin{center}
\btab{c}
\begin{picture}(60,35)(0,30)
\GgeLL{30pt}{$x_1$}
\GegLL{30pt}{$x_2$}
\GceLL{30pt}{$x_3$}
\GcgLL{30pt}{$y_1$} 
\GeeLL{30pt}{$y_2$}
\GggLL{30pt}{$y_3$}
\Gcgeg\Geggg
\Gceee\Geege
\Gcecg\Geeeg\Ggegg
\end{picture}
\\
Figure 8: $G\left(\barr{rrr}
0&1&1\\
1&1&1\\
1&1&0
\earr\right)$
\etab
\end{center}
\end{figure}

To do this we interpret the letters of $w$ differently.  If we encounter a $0$, we go around the cycle we just looped around in the same direction (clockwise or counterclockwise).  In the case of a $1$, we go around the other cycle, but keep the direction of the last cycle.  If we see a $2$, we switch cycles and direction.

For example, suppose we begin by traveling around the right-most cycle of $M$ in a clockwise direction, passing through the entries $(2,2)$, $(1,2)$, $(1,3)$, and $(2,3)$ in that order.  Now we read the first letter of $w$.  Since it is a $2$, we go around the left-most cycle of $M$ counterclockwise, passing through $(2,3)$, $(3,1)$, and $(3,2)$.  The next letter of $w$ is a $1$, so we return to the right-most cycle of $M$, but in the counterclockwise direction this time, resulting in the walk $(2,2)$, $(2,3)$, $(1,3)$, $(1,2)$.  Since the fourth letter of $w$ is a $2$, we go on to walk around the left-most cycle of $M$ in a clockwise direction.

An element of the the antichain constructed in this way is shown in Figure 7.

We have so far discussed only one manner in which antichains produced by our construction might fail to be fundamental, but in this more general setting there are a couple more subtleties.  Let us again consider the matrix
$$
\left(\barr{rrr}
0&1&1\\
1&1&1\\
1&1&0
\earr\right).
$$
If our stroll through the entries of this matrix contains a sequence like $(2,2)$, $(2,3)$, $(1,3)$, $(2,3)$, $(2,2)$, then the resulting antichain will not be fundamental.  And if our stroll does not contain infinitely many cycles, then the resulting sequence of permutations will not even form an antichain.

\section{Concluding Remarks}

Recently and independently, Albert and Atkinson~\cite{aa:sp} and Murphy~\cite{maximillian} have introduced another method for proving that closed sets of permutations are partially well-ordered.  An {\it interval\/} of $p\in S_n$ is a segment $p(i)p(i+1)\dots p(j)$, where $1\le i\le j\le n$, such that $\{p(i),\dots, p(j)\}$ forms a set of consecutive integers.  Every permutation in $S_n$ contains trivial intervals of length $1$ and $n$.  If $p$ contains no non-trivial intervals then it is said to be {\it interval-free\/} or {\it simple\/}.  For example, $35142$ is simple but $25341$ is not.

Using the full version of Higman's Theorem from \cite{higman}, it is not hard to show that an analysis of the simple permutations in a closed set can sometimes be enough to show that the set is partially well-ordered.

\begin{theorem}\label{finsimple}\cite{aa:sp,maximillian}
Let $X$ be a closed set of permutations.  If $X$ contains only finitely many simple permutations, then $X$ is partially well-ordered.
\end{theorem}

An analogous theorem for tournaments exists, and has been used in that context to show that some closed sets of permutations are partially well-ordered.  The reader is referred to Latka~\cite{latka:n5} for an example of this.

In fact, if the hypotheses of Theorem~\ref{finsimple} happen to hold, then they can be established by computer due to the following theorem, which has been proven in the special case of permutations by Albert and Atkinson~\cite{aa:sp} and Murphy~\cite{maximillian} and in the more general context of binary relational systems by Schmerl and Trotter~\cite{st}.

\begin{theorem}\cite{aa:sp,maximillian,st}
Every simple permutation of length $n>2$ contains a simple permutation of length $n-1$ or $n-2$.
\end{theorem}

Our first aim is to show that our Corollary~\ref{matrices2} is not a special case of Theorem~\ref{finsimple}.  Consider the set of permutations $D={d_1,d_2,\dots}$ given by
$$
d_k=2,4,\dots, 2k\ |\ 1,3,\dots,2k-1.
$$
(Here, as before, the vertical bar is included merely to make the permutation easier to parse, and has no mathematical meaning.)  These permutations are all simple, and $\{M_{d_1},M_{d_2},\dots\}\subset \Pr\left(\barr{r}1\\1\earr\right)$.  Hence Corollary~\ref{matrices2} (in fact, even Theorem~\ref{wk}) can be used to show that $\Pr\left(\barr{r}1\\1\earr\right)$ is partially well-ordered, whereas Theorem~\ref{finsimple} cannot draw this conclusion.

There are also situations in which Theorem~\ref{finsimple} can show that a set is partially well-ordered when Corollary~\ref{matrices2} cannot.  To describe these we need to define the {\it wreath product\/} of permutations.  Let $p$ be an $n$ permutations and let $q_1,\dots, q_n$ be permutations of any positive length.  Then $p\wr(q_1,\dots,q_n)=p_1\dots p_n$ where $p_1,\dots,p_n$ are all intervals such that
\ben
\item[(i)] $\red(p_i)=q_i$ for all $i\in[n]$, and
\item[(ii)] if $a_i$ is a term of $p_i$ for all $i\in[n]$ then $\red(a_1\dots a_n)=p$.
\een
For example, $12\wr(q_1,q_2)=q_1\oplus q_2$ and $213\wr(q_1,q_2,q_3)=(q_1\ominus q_2)\oplus q_3$.

Now let $X$ denote the largest closed set whose set of simple permutations is precisely $\{1,12,21,3142\}$.  Such a set exists because the union of any number of closed sets is again a closed set.  This set is partially well-ordered by Theorem~\ref{finsimple}, but it contains the infinite sequence of permutations
\bea
z_1&=&3,1,4,2,\\
z_2&=&11, 9, 12, 10\ |\ 3, 1, 4, 2\ |\ 15, 13, 16, 14\ |\ 7, 5, 8, 6,\\
&\dots&\\
z_k&=&z_{k-1}\wr(3142,3142,\dots,3142).
\eea
It is perhaps easiest to see the pattern in this sequence by looking at the matrices, for example,
$$
M_{z_2}=
\left(
\barr{rrrrrrrrrrrrrrrr}
&&&&&&&&&&1&&&&&\\
&&&&&&&&1&&&&&&&\\
&&&&&&&&&&&1&&&&\\
&&&&&&&&&1&&&&&&\\
&&1&&&&&&&&&&&&&\\
1&&&&&&&&&&&&&&&\\
&&&1&&&&&&&&&&&&\\
&1&&&&&&&&&&&&&&\\
&&&&&&&&&&&&&&1&\\
&&&&&&&&&&&&1&&&\\
&&&&&&&&&&&&&&&1\\
&&&&&&&&&&&&&1&&\\
&&&&&&1&&&&&&&&&\\
&&&&1&&&&&&&&&&&\\
&&&&&&&1&&&&&&&&\\
&&&&&1&&&&&&&&&&
\earr
\right).
$$
Corollary~\ref{matrices2} cannot conclude that $X$ is partially well-ordered because there is no finite $\zpm$ matrix $M$ for which $M_{z_k}$ lies in the strong completion of $\Pr(M)$ for all $k$.

To complete our examples, consider the set of permutations gotten by reducing the first $2k$ terms of the infinite sequence
$$
4,1,6,3,8,5,10,7,\dots,
$$
referred to in \cite{maximillian} as the increasing oscillating sequence.  Clearly the closure of the resulting set is partially well-ordered, since it is the closure of a chain.  However, every one of these permutations is simple, so Theorem~\ref{finsimple} cannot reach this conclusion, and it is easy to see that Corollary~\ref{matrices2} is also of no help.

Given these examples, it seems natural to ask if there a common generalization of these two techniques.  We do not have an answer for this.

\bigskip

\noindent{\it Acknowledgment.\/}  The authors met at the 2003 Conference on Permutation Patterns.  Both of their visits were partially supported by the New Zealand Institute of Mathematics and its Applications.  In addition, Vince Vatter thanks Doron Zeilberger for support.

\begin{\bib}{99}

\bibitem{aa:sp} M. H. Albert and M. D. Atkinson, Simple permutations and pattern restricted permutations, in preparation.

\bibitem{a:rp} M. D. Atkinson, Restricted permutations, {\it Discrete Math.\/} {\bf 195} (1999), 27-38.

\bibitem{amr:pwocsop} M. D. Atkinson, M. M. Murphy, and M. Ru\v{s}kuc, Partially well-ordered closed sets of permutations, {\it Order\/} {\bf 19} (2002), 101-113.

\bibitem{as:rpwp} M. D. Atkinson and T. Stitt, Restricted permutations and the wreath product, {\it Discrete Math.\/} {\bf 259} (2002), 19-36.

\bibitem{cl:minimal} G. L. Cherlin and B. J. Latka, Minimal antichains in well-founded quasi-orders with an application to tournaments, {\it J. Combin. Theory Ser. B\/} {\bf 80} (2000), 258--276.

\bibitem{fh:ds} Z. F\H{u}redi and P. Hajnal, Davenport-Schinzel theory of matrices, {\it Discrete Math.\/} {\bf 103} (1992), 233-251.

\bibitem{gustedt} J. Gustedt, ``Algorithmic Aspects of Ordered Structures,'' Ph.D. dissertation, Technische Universit\"at Berlin, 1992.

\bibitem{higman} G. Higman, Ordering by divisibility in abstract algebra, {\it Proc. London Math. Soc.\/} {\bf 2} (1952), 326-336.

\bibitem{latka:n5} B. J. Latka, Tournaments that omit $N_5$ are well-quasi-ordered, preprint.

\bibitem{laver} R. Laver, Well-quasi-orderings and sets of finite sequences, {\it Math. Proc. Camb. Philos. Soc.\/} {\bf 79} (1976), 1-10.

\bibitem{lothaire1} M. Lothaire, ``Combinatorics on Words,'' Encyclopedia of Math., Vol.\ 17, Addison Wesley, Reading, MA, 1983.

\bibitem{maximillian} M. M. Murphy, Ph.D. dissertation, University of St. Andrews, 2002.

\bibitem{p:cpwdeqpsapq} V. R. Pratt, Computing permutations with double-ended queues, parallel stacks and parallel queues, {\it Proc. ACM Symp. Theory of Computing\/} {\bf 5} (1973), 268-277.

\bibitem{st} J. Schmerl and W. Trotter, Critically indecomposable partially ordered sets, graphs, tournaments and other binary relational structures, {\it Discrete Math.\/} {\bf 113} (1993), 191-205.

\bibitem{ss:rp} R. Simion and F. W. Schmidt, Restricted permutations, {\it European J. Combin.\/} {\bf 6} (1985), 383--406.

\bibitem{sb:aiaop} D. Spielman, M. B\'ona, An infinite antichain of permutations, {\it Electronic J. Combinatorics\/} {\bf 7} (2000), \#N2.

\bibitem{t:sunoqas} R. E. Tarjan, Sorting using networks of queues and stacks, {\it J. of the ACM\/} {\bf 19} (1972), 341-346.

\bibitem{w:pip} H. Wilf, The patterns of permutations, {\it Discrete Math.\/} {\bf 257} (2002), 575-583.

\end{\bib}

\end{document}